\documentclass[aos,MSNbibl,dvips]{arximspdf}

\doi{10.1214/14-AOS1254} 
\volume{42}
\issue{6}
\pubyear{2014}
\firstpage{2340}
\lastpage{2381}
\docsubty{FLA}

\makeatletter

\newcommand{\eqref}[1]{(\ref{#1})}
\newcommand{\ep}{\varepsilon}
\newcommand{\rrvert}{\vert}
\newcommand{\rrVert}{\Vert}
\newcommand{\llvert}{\vert}
\newcommand{\llVert}{\Vert}
\renewcommand{\mid}{|}
\newtheorem{teo}{Theorem}[section]
\newtheorem{lmm}[teo]{Lemma}
\newtheorem{cor}[teo]{Corollary}
\newtheorem{prop}[teo]{Proposition}
\newcommand{\mg}{\mathcal{G}}
\newcommand{\cov}{\operatorname{Cov}}
\newcommand{\ee}{\mathbb{E}}
\newcommand{\mf}{\mathcal{F}}
\newcommand{\pp}{\mathbb{P}}
\newcommand{\rr}{\mathbb{R}}
\newcommand{\var}{\operatorname{Var}}
\newcommand{\hmm}{\hat{\mu}}
\newcommand{\tm}{\tilde{\mu}}
\makeatother

\begin{document}
\begin{frontmatter}

\title{A new perspective on least squares under convex~constraint\thanksref{T1}}
\runtitle{Least squares under convex constraint}

\begin{aug}
\author[A]{\fnms{Sourav}~\snm{Chatterjee}\corref{}\ead[label=e1]{souravc@stanford.edu}}
\runauthor{S. Chatterjee}
\affiliation{Stanford University}
\address[A]{Department of Statistics\\
Stanford University\\
Sequoia Hall, 390 Serra Mall\\
Stanford, California 94305\\
USA\\
\printead{e1}} 
\end{aug}
\thankstext{T1}{Supported in part by NSF Grant DMS-10-05312.}

\received{\smonth{3} \syear{2014}}
\revised{\smonth{6} \syear{2014}}

%
\begin{abstract}
Consider the problem of estimating the mean of a Gaussian random vector
when the mean vector is assumed to be in a given convex set. The most
natural solution is to take the Euclidean projection of the data vector
on to this convex set; in other words, performing ``least squares under
a convex constraint.'' Many problems in modern statistics and statistical
signal processing theory are special cases of this general situation.
Examples include the lasso and other high-dimensional regression techniques,
function estimation problems, matrix estimation and completion,
shape-restricted regression, constrained denoising, linear inverse problems,
etc. This paper presents three general results about this problem, namely,
(a)~an exact computation of the main term in the estimation error by relating
it to expected maxima of Gaussian processes (existing results only give
upper bounds),
(b)~a theorem showing that the least squares estimator is always
admissible up to a
universal \mbox{constant} in any problem of the above kind and (c) a
counterexample showing
that least squares estimator may not always be minimax rate-optimal.
The result from
part~(a) is then used to compute the error of the least squares
estimator in two
examples of contemporary interest.
\end{abstract}

%
\end{frontmatter}

\section{Theory}\label{theory}
\subsection{The setup}
Throughout this manuscript, $Z= (Z_1,\ldots, Z_n)$ denotes an
\mbox{$n$-}dimensional standard Gaussian random vector. Let $\mu=(\mu
_1,\ldots, \mu_n)\in\rr^n$ be a point in $\rr^n$, and let $Y=
Z+\mu$. We are interested in estimating $\mu$ from the data vector
$Y$. If nothing more is known, the vector $Y$ itself is the maximum
likelihood estimate of $\mu$.

Suppose now that $\mu$ is known to belong to a closed convex set
$K\subseteq\rr^n$. Let $P_K$ denote the Euclidean projection on to
$K$. That is, for a vector $x\in\rr^n$, $P_K(x)$ is the point in $K$
that is closest to $x$ in the Euclidean distance. It is a standard fact
about closed convex sets (see Lemma~\ref{proj} in Section~\ref
{proofs}) that $P_K$ is a well-defined map. Under the assumption that
$\mu\in K$, the maximum likelihood estimate of $\mu$ in the Gaussian
model is $\hmm:= P_K(Y)$. We will refer to $\hmm$ as the least
squares estimator (LSE) of $\mu$ under the convex constraint~$K$.
As mentioned in the abstract, many problems in modern statistics are
special cases of this general setup, including the lasso and other
high-dimensional regression techniques, function estimation problems,
matrix estimation and completion, shape-restricted regression, etc.

Let $\llVert  x\rrVert  $ denote the Euclidean norm of a vector $x\in\rr^n$. Our
main goal is to understand the magnitude of the estimation error $\llVert  \hmm- \mu\rrVert $. The standard approach to computing upper bounds on the
expected squared value of this error (the ``risk'') is via empirical
process theory and related entropy computations. As a consequence of
path-breaking contributions from a number of authors over a period of
more than thirty years, including Birg\'e \cite{birge83}, Tsirelson
\cite{tsirel1,tsirel2,tsirel3}, Pollard \cite{pollard}, van de Geer
\cite{vdg87,vdg90,vdg93}, Birg\'e and Massart \cite{birgemassart},
van der Vaart and Wellner \cite{vdvwellner} and many others, we now
have a fairly good idea about how to convert results for expected
maxima of empirical processes to upper bounds on estimation errors in
problems of the above type, especially in the context of regression. To
know more about this important branch of theoretical statistics and
machine learning, see the monographs of B\"uhlmann and van de Geer
\cite{buhlmannvdg}, Massart \cite{massart}, van de Geer \cite
{vdgbook} and van der Vaart and Wellner \cite{vdvwellner}.

In a different direction, this general problem has recently gained
prominence in the statistical signal processing literature. The least
squares problem outlined above is essentially equivalent to the problem
of constrained denoising in signal processing. It is also intimately
connected with the so- called linear inverse problems. The history of
this line of investigation, although of relatively recent origin, is
already quite formidable. Important papers include those of Rudelson
and Vershynin \cite{rv08}, Stojnic \cite{stojnic09}, Oymak and
Hassibi \cite{oh10,oh13}, Chandrasekaran et~al.~\cite
{chandrasekaran}, Chandrasekaran and Jordan \cite{cj13} and Amelunxen
et~al.~\cite{amelunxen}. For some interesting recent developments, see
McCoy and Tropp \cite{mt13a,mt13b}, Foygel and Mackey \cite{fm} and
Tropp \cite{tropp14}.

In the signal processing context, the expected squared error $\ee\llVert \hmm-\mu\rrVert  ^2$ is closely related to the concept of ``statistical
dimension'' recently introduced in Amelunxen et~al.~\cite{amelunxen}.
It is also related to the older existing notions of ``Gaussian width''
in probability, functional analysis and signal processing,
``mean width'' in convex geometry and ``Gaussian complexity'' in
learning theory and machine learning.

\subsection{Main result}\label{mainestsec}
One limitation of the theory based on empirical processes in its
current form is that it only gives upper bounds on the error. There are
some lower bounds
``in spirit,'' in the form of necessary and sufficient conditions for
consistency (e.g., in Tsirelson \cite{tsirel1} and van de Geer and
Wegkamp \cite{vdgwegkamp}) but the lower bounds are not explicit. The
first main result of this manuscript, presented below, shows that if
one looks at expected maxima of certain Gaussian processes (instead of
upper bounds on these maxima) then one can get an approximation for the
actual error instead of just an upper bound. Not only that, the theorem
also shows that the error $\llVert  \hmm-\mu\rrVert  $ is typically concentrated
around its expected value.

Let $x\cdot y$ denote the usual inner product on $\rr^n$ and let $K$
be any nonempty closed convex set. For any $\mu\in\rr^n$ and any
$t\ge0$, let
\[
f_\mu(t):= \ee \Bigl(\sup_{\nu\in K\dvtx   \llVert  \nu-\mu\rrVert  \le t} Z\cdot(\nu-\mu)
\Bigr) - \frac{t^2}{2},
\]
where $Z$ is an $n$-dimensional standard Gaussian random vector. If
$\mu\notin K$, then there is no $\nu\in K$ satisfying $\llVert  \mu-\nu\rrVert
\le t$ if $t$ is strictly less than the distance of $\mu$ from $K$. In
that case, define $f_\mu(t)$ to be $-\infty$, following the standard
convention that the supremum of an empty set is $-\infty$.

Let $t_\mu$ be the point in $[0,\infty)$ where $f_\mu$ attains its
maximum. We will show below that $t_\mu$ exists and is unique. Recall
that $P_K$ denotes the projection on to $K$, and that
\[
\hmm:= P_K(Z+\mu)
\]
is the least squares estimate of $\mu$ based on the data vector $Z+\mu
$. The following theorem shows that irrespective of the dimension $n$
and the convex set $K$, it is always true that
\[
\llVert \hmm- \mu\rrVert = t_\mu+ O\bigl(\max\{\sqrt{t_\mu},
1\}\bigr).
\]
In particular, if $t_\mu$ is large, then the random quantity $\llVert  \hmm
-\mu\rrVert  $ is concentrated around the nonrandom value $t_\mu$. 

%
\begin{teo}\label{mainest}
Let $K$, $\mu$, $\hmm$, $f_\mu$ and $t_\mu$ be as above. Let $t_c:= \inf_{\nu\in K} \llVert  \nu-\mu\rrVert  $.
Then $f_\mu(t)$ is equal to $-\infty$ when $t<t_c$, is a finite and
strictly concave function of $t$ when $t\in[t_c,\infty)$, and decays
to $-\infty$ as $t\rightarrow\infty$. Consequently, $t_\mu$ exists
and is unique. Moreover, for any $x\ge0$,
\[
\pp \bigl( \bigl\llvert \llVert \hmm-\mu\rrVert - t_\mu \bigr\rrvert
\ge x\sqrt{t_\mu } \bigr)\le3\exp \biggl(-\frac{x^4}{32(1+{x}/{\sqrt{t_\mu
}})^2} \biggr).
\]
\end{teo}
Note that $\mu$ is not required to be in $K$ in this theorem. The tail
bound is valid even if $\mu$ is a point lying outside $K$.

The above theorem can potentially give rise to many corollaries. One
basic corollary, presented below, gives estimates for the expected
squared error of $\hmm$. Although Theorem~\ref{mainest} contains a
lot more information than this corollary, expected squared errors are
culturally important.

%
\begin{cor}\label{maincor}
Let all notation be as in Theorem~\ref{mainest}. Then there is a
universal constant $C$ such that if $t_\mu\ge1$, then
\[
t_\mu^2 - C t_\mu^{3/2}\le\ee
\llVert \hmm- \mu\rrVert ^2 \le t_\mu^2 + C
t_\mu^{3/2},
\]
and if $t_\mu< 1$, then
\[
\ee\llVert \hmm-\mu\rrVert ^2 \le C.
\]
\end{cor}
It may be illuminating to see an example at this point. Consider the
simplest possible example, namely, that $K$ is a $p$-dimensional
subspace of $\rr^n$, where $p\le n$. This is nothing but the linear
regression setup, assuming that $\mu= X\beta$, where $X$ is an
$n\times p$ matrix of full rank and $\beta\in\rr^p$ is arbitrary.

Since $K$ is a subspace, $Z\cdot x = P_K(Z) \cdot x$ for any $x\in K$.
Moreover, $P_K(Z)$ is a standard Gaussian random vector in $K$. A
simple application of the rotational invariance of $Z$ shows that we
may assume, without loss of generality, that $K$ is simply a copy of
$\rr^p$ contained in $\rr^n$. Combining these observations, we see
that for any $\mu\in K$ and $t\ge0$,
\begin{eqnarray*}
f_\mu(t) &=& \ee \Bigl(\sup_{x\in\rr^p,   \llVert  x\rrVert  \le t} W\cdot x \Bigr)
- \frac{t^2}{2},
\end{eqnarray*}
where $W$ is a $p$-dimensional standard Gaussian random vector. The
above expression can be exactly evaluated, to give
\[
f_\mu(t) = \ee\bigl(t\llVert W\rrVert \bigr) -\frac{t^2}{2}.
\]
Clearly, $f_\mu$ is maximized at
\[
t_\mu= \ee\llVert W\rrVert = \sqrt{p} + O(1),
\]
where $O(1)$ denotes a quantity that may be bounded by a constant that
does not depend on $p$ or $n$. By Theorem~\ref{mainest}, this shows
that when $K$ is a $p$-dimensional subspace of $\rr^n$, then with high
probability,
\[
\llVert \hmm- \mu\rrVert = \sqrt{p} + O\bigl(p^{1/4}\bigr).
\]
Of course, this result may be derived by other means. It is included
here only to serve as a simple illustration.

The above example is, in some sense, exceptionally simple. In general,
it will be very difficult to compute $t_\mu$ exactly, since we have
only limited tools at our disposal to compute expected maxima of
high-dimensional Gaussian processes. However, the strict concavity of
the function $f_\mu$ gives an easy way to calculate upper and lower
bounds on $t_\mu$ (and hence, upper and lower bounds on the estimation
error $\llVert  \hmm-\mu\rrVert  $) by calculating bounds on $f_\mu$ at a small
number of points.

%
\begin{prop}\label{simple}
If $0\le r_1< r_2$ are such that $f_\mu(r_1)\le f_\mu(r_2)$, then
$t_\mu\ge r_1$. On the other hand, if $f_\mu(r_1)\ge f_\mu(r_2)$,
then $t_\mu\le r_2$. In particular, if $\mu\in K$ and $r>0$ is such
that $f_\mu(r)\le0$, then $t_\mu\le r$.
\end{prop}
In Section~\ref{applications}, we will see applications of this
proposition in computing matching upper and lower bounds for estimation
errors in two nontrivial problems.

\subsection{The LSE is admissible up to a universal constant}
Does the constrained least squares estimator $\hmm$ enjoy any kind of
general optimality property that holds for any $K$? A priori, this may
sound like a hopeless question due to the completely arbitrary nature
of the convex set $K$. One may hope that the LSE is minimax optimal in
some asymptotic sense, but as we will show later, this is not the case.
Fortunately, it turns out that $\hmm$ indeed enjoys a certain other
kind of universal optimality property, as shown in Theorem~\ref
{admissibility} below. From a purely mathematical point of view, this
is the deepest result of this paper.

The famous Stein paradox~\cite{stein} shows that the least squares
estimate $\hmm$ is inadmissible under square loss when $K=\rr^n$.
Stein's example gave birth to the flourishing field of shrinkage
estimates. The following theorem shows that although the LSE $\hmm$
may be inadmissible, it is always
``admissible up to a universal constant,'' whatever be the set $K$. In
particular, shrinkage---or any other clever idea---cannot improve
the risk beyond a universal constant factor everywhere on the parameter
space. 

%
\begin{teo}\label{admissibility}
There is a universal constant $C>0$ such that the following is true.
Take any $n$ and any nonempty closed convex set $K\subseteq\rr^n$.
Let \mbox{$g\dvtx \rr^n \rightarrow\rr^n$} be any Borel measurable map, and for
each $\mu\in\rr^n$ define the estimate $\tilde{\mu}:= g(Z+\mu)$,
where $Z$ is a standard Gaussian random vector. Let $\hmm$ be the
least squares estimate $P_K(Z+\mu)$, as in Theorem~\ref{mainest}.
Then there exists $\mu\in K$ such that $\ee\llVert  \tm- \mu\rrVert  ^2 \ge C
\ee\llVert  \hmm- \mu\rrVert  ^2$.
\end{teo}
Again, it may be a good idea to understand the impact of Theorem~\ref
{admissibility} through an example. Consider the problem of $\ell
^1$-penalized regression with $p$ covariates, where $p$ may be bigger
than $n$. Here, $K$ is the set of all $\mu$ of the form $X\beta$,
where $X$ is a given $n\times p$ matrix and $\beta$ is a point in $\rr
^p$ with $\ell^1$ norm bounded by some prespecified constant~$L$. The
convex-constrained least squares estimate in this problem is the same
as the lasso estimate of Tibshirani \cite{tibs96} in its primal form.
One may consider various other procedures for computing estimates of
$\beta$ in this problem. Theorem~\ref{admissibility} says that no
matter what procedure one considers, there is always some $\beta$ with
$\ell^1$ norm $\le L$ where the prediction error of the new procedure
is at least as big as the prediction error of the lasso, multiplied by
a universal constant.

It is interesting to figure out the optimal value of the universal
constant in Theorem~\ref{admissibility}. Note that by the Stein
paradox, the largest possible value is strictly less than $1$.

\subsection{The LSE may not be minimax rate-optimal}\label{countersec1}
Theorem~\ref{admissibility} shows that there is always some region of
the parameter space where the least squares estimate $\hmm$ does not
perform too badly in comparison to any given competitor. This
immediately raises the question as to whether the same is true about
the maximum risk: is the maximum risk of the least squares estimate
always within a universal constant multiple of the minimax risk? (Here the
``risk'' of an estimator $\tilde{\mu}$ under square loss is defined,
as usual, to be $\ee\llVert  \tilde{\mu}-\mu\rrVert  ^2$.) Surprisingly, the
answer turns out to be negative, as shown by the following counterexample.

Take any $n$. Define a closed convex set $K\subseteq\rr^n$ as
follows: take any $\alpha\in[0,1]$, $\theta_1,\ldots,\theta_n\in
[-1,1]$, and let
\[
\mu_i:= \alpha n^{-1/4} + \alpha\theta_i
n^{-1/2},\qquad i=1,\ldots, n.
\]
Let $K$ be the set of all $\mu= (\mu_1,\ldots,\mu_n)$ obtained as above.

%
\begin{prop}\label{counter1}
The set $K$ defined above is closed and convex. As before, let $\hmm=
P_K(Z+\mu)$ be the least squares estimate of $\mu\in K$ obtained by
projecting the data vector $Y=Z+\mu$ on to $K$. Let $\tilde{\mu}$ be
the estimate whose coordinates are all equal to the average of the
coordinates of $Y$. Then, under square loss, the maximum risk of $\hmm
$ is bounded below by $C_1n^{1/2}$ whereas the maximum risk of $\tilde
{\mu}$ is bounded above by $C_2$, where $C_1$ and $C_2$ are positive
constants that do not depend on~$n$.
\end{prop}
It is interesting to understand whether this example is a pathological
exception, or if there is a general rule that dictates whether the LSE
is minimax rate-optimal or not in a given problem. Theorem~\ref
{admissibility} gives a sufficient condition for minimax
rate-optimality, namely, that the risk is of the same order everywhere
on $K$. This is expressed quantitatively in the following proposition.
But this condition may be difficult to verify in examples.

%
\begin{prop}\label{adcor}
Let $\tm$ and $\hmm$ be as in Theorem~\ref{admissibility}. For each
$\mu\in K$, let $R_1(\mu)$ be the risk of $\hmm$ at $\mu$ and
$R_2(\mu)$ be the risk of $\tm$ at $\mu$. Then
\[
\sup_{\mu\in K} R_2(\mu) \ge CL \sup
_{\mu\in K} R_1(\mu),
\]
where $C$ is the universal constant from Theorem~\ref{admissibility} and
\[
L:= \frac{ \inf_{\mu\in K} R_1(\mu)}{ \sup_{\mu\in K}R_1(\mu
)}.
\]
\end{prop}
The counterexample given in Proposition~\ref{counter1} also raises the
question as to whether there is a general estimator that is guaranteed
to be minimax up to a universal constant.

\section{Examples}\label{applications}
This section contains two nontrivial applications of Theorem~\ref
{mainest}, to supplement the easy example worked out in Section~\ref
{mainestsec}. We only present the results here. The details are worked
out in Section~\ref{proofs}. 

\subsection{Lasso with nonsingular design}
Let $p\ge1$ and $n\ge2$ be two integers, and let $X$ be a given
$n\times p$ matrix with real entries. Let $L$ be a positive real
number, and let
%
\begin{equation}
\label{lassok0} K_0:= \bigl\{\beta\in\rr^p\dvtx  \llvert
\beta\rrvert _1 \le L\bigr\},
\end{equation}
where $\llvert  \beta\rrvert  _1$ stands for the $\ell^1$ norm of $\beta$, that is,
the sum of the absolute values of the components of $\beta$. Let
%
\begin{equation}
\label{lassok} K:= \{X\beta\dvtx  \beta\in K_0\}.
\end{equation}
The least squares estimator for the convex constraint $K$ is nothing
but the lasso estimator in its primal form as defined by Tibshirani
\cite{tibs96}. The number $L$ is called the
``penalty parameter.''

The theoretical properties of the lasso and related procedures have
been extensively studied over the last ten years, notably by Donoho and
coauthors \cite
{donoho,donoho04,donohoelad02,donohoetal,donohohuo02,donohojohnstone94,donohoetal95},
Knight and Fu \cite{fu},
Zou \cite{zou}, Wainwright \cite{wainwright}, Cand\`es and Tao \cite
{candestao05,candestao07}, Meinshausen and B\"uhlmann \cite{meins1},
Meinshausen and Yu \cite{meins2}, Koltchinskii \cite{kol09}, Wang and
Leng \cite{wangleng}, Zhao and Yu \cite{zhaoyu}, Bunea et~al.~\cite
{bunea}, van de Geer \cite{vdg08}, Greenshtein and Ritov \cite{gr},
Bickel et~al.~\cite{bickeletal09}, Bartlett et~al.~\cite
{bartlettetal12}, Rigollet and Tsybakov \cite{rigollet}, Oymak et~al.~\cite{oymak} and many others. For a more complete set of
references and a clear exposition of the results and techniques, see
the wonderful recent monograph of B\"uhlmann and van de Geer \cite
{buhlmannvdg}.
The investigators have tried to understand a number of different kinds
of consistency for the lasso estimator. The expected squared error $\ee
\llVert  \hmm-\mu\rrVert  ^2$ translates into what is known as the
``squared prediction error'' in the lasso literature. Among the papers
cited above, the ones dealing mainly with the behavior of the
prediction error are \cite{bunea,vdg08,gr}. If the prediction error
vanishes on an appropriate scale, the lasso procedure is called
``risk consistent.''

Risk consistency does not require too many assumptions~\cite
{chat,buhlmannvdg,gr}, but the available bounds on the expected squared
prediction error are solely upper bounds. Matching lower bounds are not
known in any case. In particular, it is not known how the error depends
on the choice of the penalty parameter $L$. Practitioners believe from
experience that choosing the penalty parameter correctly is of crucial
importance, and this is usually done using cross-validation of some
sort, for example, in Tibshirani \cite{tibs96,tibs11}, Greenshtein and
Ritov \cite{gr}, Hastie et~al.~\cite{hastie}, Efron et~al.~\cite
{efron}, and van de Geer and Lederer \cite{vdgled}, although some
other techniques have also been proposed, for example, in Tibshirani
and Taylor \cite{tt} and Zou et~al.~\cite{zouetal}. For some nascent
theoretical progress on cross-validation for the lasso and further
references, see Homrighausen and McDonald~\cite{homrig}.

The following theorem demonstrates, for the first time, the critical
importance of choosing the correct penalty parameter value. If the
penalty parameter $L$ is chosen to be equal to $\llvert  \beta\rrvert  _1$, then the
prediction error is vastly smaller than if the two quantities are
unequal. Although the theorem is restricted to the case of nonsingular
design matrices, we may expect the phenomenon to hold in greater generality.

%
\begin{teo}\label{lassothm}
Take any $L>0$ and let $K$ be defined as in~(\ref{lassok}). Let
$\Sigma:= X^TX/n$, and let $a$ and $b$ be the smallest and largest
eigenvalues of $\Sigma$. Assume that $a>0$, and that all the diagonal
entries of $\Sigma$ are equal to $1$. Take any $\beta\in\rr^p$ and
let $\mu:=X\beta$. Let $s$ be the number of nonzero entries of $\beta
$. Let
\[
\delta:= L-\llvert \beta\rrvert _1
\]
and $r:= p/n$. Let $t_\mu$ be as in Theorem~\ref{mainest}, for the
set $K$ defined in (\ref{lassok}). If $\delta> 0$, then given any
$\varepsilon> 0$ there is a constant $C_1$ depending only on $\delta$,
$\varepsilon$, $a$, $b$, $s$, $r$ and~$L$ such that whenever $n> C_1$, we have
\[
n^{1/4-\varepsilon} \le t_\mu\le n^{1/4+\varepsilon}.
\]
If $\delta= 0$, then there is a constant $C_2$ depending only on $a$,
$b$, $s$, $r$ and $L$ such that
\[
t_\mu\le C_2 \sqrt{\log n}.
\]
Finally, if $\delta< 0$, there are positive constants $C_3$ and $C_4$
depending only on $\delta$, $a$, $b$, $s$, $r$ and $L$ such that
\[
C_3 \sqrt{n}\le t_\mu\le C_4 \sqrt{n}.
\]
\end{teo}
The reader may easily check the implications of the above bounds on the
prediction error by looking back at Theorem~\ref{mainest}. In
particular,\vspace*{1pt} they show that the squared prediction error $\ee\llVert  X
\hat{\beta}-X\beta\rrVert  ^2$ equals $n^{1/2+o(1)}$ if the penalty
parameter $L$ is greater than $\llvert  \beta\rrvert  _1$, is of order $n$ if $L$ is
less than $\llvert  \beta\rrvert  _1$, and is bounded above by some constant multiple
of $\log n$ if the penalty parameter is chosen correctly, to be equal
to $\llvert  \beta\rrvert  _1$. Therefore, it is very important that $L$ is chosen
correctly when implementing the lasso procedure. However, there is a
caveat: Theorem~\ref{lassothm} does not prove anything in the case
where $L$ is chosen using the data. It only shows the importance of
choosing the correct value of the penalty parameter, besides being the
first result that establishes a lower bound on the lasso error. To
prove an analogous result for the case where $L$ is chosen using the
data requires further work.

\subsection{Isotonic regression}
Define the convex set
%
\begin{equation}
\label{isok} K:= \bigl\{(\mu_1,\ldots,\mu_n)\in
\rr^n\dvtx  \mu_1\le\mu_2\le\cdots \le
\mu_n\bigr\}.
\end{equation}
The least squares problem for this convex constraint, popularly known
as ``isotonic regression'' or ``monotone regression,'' has a long
history in the statistics literature, possibly beginning in Ayer et~al.~\cite{ayeretal} and Grenander \cite{grenander}. The LSE is easily
computed using the so-called ``pool adjusted violators algorithm'' (see
Robertson et~al.~\cite{robertson}, Chapter~1).

There is substantial literature on the properties of individual $\hmm
_i$, as $i/n$ is fixed and $n$ goes to infinity, with some appropriate
limiting behavior assumed for the mean vector $\mu$. Some notable
papers on such local errors are those of Prakasa Rao \cite{prakasa},
Brunk~\cite{brunk}, Groeneboom and Pyke~\cite{gp}, Durot~\cite
{durot}, Carolan and Dykstra~\cite{carolan}, Cator~\cite{cator} and
Jankowski \cite{jw}. The global error $\llVert  \hmm-\mu\rrVert  $ has
also received considerable attention, notably in van de Geer \cite
{vdg90,vdg93}, Donoho \cite{donoho}, Birg\'e and Massart \cite
{birgemassart}, Wang \cite{wang}, Meyer and Woodroofe \cite{mw},
Zhang \cite{zhang} and Chatterjee, Guntoboyina and Sen \cite{cgs}.

It is now generally understood that if the $\mu_i$'s are ``strictly
increasing'' in some limiting sense, then $\hmm_i-\mu_i$ is typically
of order $n^{-1/3}$, whereas the error is smaller if the $\mu_i$'s
have ``flat stretches'' \cite{cgs}. Therefore, it is natural to expect
that in the strictly increasing case, $\llVert  \hmm-\mu\rrVert  $ should be of
order $n^{1/6}$. Using Theorem~\ref{mainest}, it turns out that we may
not only get finite sample upper and lower bounds for the global risk
$\ee\llVert  \hmm-\mu\rrVert  ^2$, but also show that $\llVert  \hmm-\mu\rrVert  $ is
concentrated around its mean value; that is, there is some constant
$C(\mu)$ depending on $\mu$ such that with high probability,
%
\begin{equation}
\label{isobds} \llVert \hmm-\mu\rrVert = C(\mu)n^{1/6} + O
\bigl(n^{1/12}\bigr).
\end{equation}
The following theorem makes this precise.

%
\begin{teo}\label{convthm}
Let $K$ be the convex set defined in (\ref{isok}). Take any $\mu\in
K$ and let $\hmm= P_K(Z+\mu)$ be the LSE of $\mu$ obtained from the
data vector $Z+\mu$. Let
\begin{eqnarray*}
D &:=& \max\{\mu_n-\mu_1, 1\},
\\
A &:=& \min_{1\le i\le n-1} n(\mu_{i+1}-\mu_i),
\\
B &:=& \max_{1\le i\le n-1} n(\mu_{i+1}-\mu_i).
\end{eqnarray*}
Let $t_\mu$ be as in Theorem~\ref{mainest}, for the set $K$ defined
in (\ref{isok}). Then
\[
\frac{C_1 A^{8/3} n^{1/6}}{B^{4/3}D}\le t_\mu\le C_2 D^{1/3}
n^{1/6},
\]
where $C_1$ and $C_2$ are positive universal constants.
\end{teo}
The reader may easily check the consequences of the above bounds on
$t_\mu$ by looking back at Theorem~\ref{mainest} and Corollary~\ref
{maincor}, and in particular, that it proves (\ref{isobds}) when $D$,
$A$ and $B$ are all of constant order.

Just to be clear, the upper bound on the expected mean-squared error
that we get from Theorem~\ref{convthm} can be derived from existing
results such as those in Zhang \cite{zhang}. The new contribution of
Theorem~\ref{convthm} is the lower bound, and also the conclusion (in
combination with Theorem~\ref{mainest}) that the squared error
concentrates around its expected value.

\section{Proof sketches}
Since the proofs of the main results (Theorems~\ref{mainest} and~\ref{admissibility}) are somewhat technical, I will try to give a readable
sketch of the main ideas in this section. The details are given in
Section~\ref{proofs}.

The proof of Theorem~\ref{mainest} goes roughly as follows. Define a
random function
\[
F_\mu(t) = \sup_{\nu\in K\dvtx   \llVert  \nu-\mu\rrVert  \le t} Z\cdot(\nu-\mu ) -
\frac{t^2}{2},
\]
so that $f_\mu(t) = \ee(F_\mu(t))$. Using the convexity of $K$,
prove that $F_\mu$ and $f_\mu$ are both strictly concave functions.
Let $t^*$ be the unique point at which $F_\mu$ is maximized. Again,
use convexity of $K$ and some algebraic manipulation to prove the key identity
\[
\llVert \hmm-\mu\rrVert = t^*.
\]
Note that this is a purely deterministic identity, having nothing to do
with the modeling assumptions.

Next, using the concentration of Gaussian maxima, show that $F_\mu
(t)$, although random, is close to $f_\mu(t)$ with high probability.
Since $F_\mu$ and $f_\mu$ are two strictly concave functions that are
close to each other with high probability, their points of maxima must
also be close. That is, $t_\mu\approx t^* = \llVert  \hmm-\mu\rrVert  $ with high
probability.

The proof of Theorem~\ref{admissibility} is more complex and it is
quite hard to present the ideas in a nutshell. Still, a high level
overview of the main steps may be given as follows.

Throughout this proof sketch, ``constant'' will mean ``positive
universal constant.'' Take any $\mu^* \in K$ and let $B_0$ be a ball
of radius $C_1t_{\mu^*}$ around $\mu^*$, where $C_1$ is a small
constant that will be chosen later. Let $\rho$ be a probability
measure on $B_0$, also to be chosen later. Let $g(Z+\mu)$ be any
estimator of $\mu$. Suppose that we are able to prove
%
\begin{equation}
\label{s1} \int_{B_0} \ee\bigl\llVert \mu-g(Z+\mu)\bigr
\rrVert ^2 \,d\rho(\mu)\ge C_2 t_{\mu^*}^2
\end{equation}
for some constant $C_2$. Upon solving some technical hurdles, it can be
shown that the dependence of $t_\mu$ on $\mu$ is smooth enough to
guarantee that if the constant $C_1$ is chosen small enough, then there
is a constant $C_3$ such that $t_{\mu^*} \ge C_3 t_\mu$
for all $\mu\in B_0$. Combined with (\ref{s1}) and Corollary~\ref
{maincor}, this implies the existence of $\mu_0\in B_0$ such that
\[
\ee\bigl\llVert \mu_0-g(Z+\mu_0)\bigr\rrVert
^2\ge C_2t_{\mu^*}^2 \ge
C_2 C_3^2 t_{\mu
_0}^2
\ge C_4 \ee\bigl\llVert \mu_0-P_K(Z+
\mu_0)\bigr\rrVert ^2
\]
for some constant $C_4$, completing the proof.

The main challenge, therefore, is to show (\ref{s1}). We now make a
specific choice of $\rho$. Let $\rho$ be the probability measure of
the point $\nu^*$ that maximizes $Z\cdot(\mu-\mu^*)$ among all $\mu
\in B_0$. Let $Z'$ be an independent copy of $Z$, and let $Y' = Z'+\nu
^*$. Then observe that the expression on the left-hand side of (\ref
{s1}) is nothing but $\ee\llVert  \nu^* - g(Y')\rrVert  ^2$.

The main trick now is the following. Let $\nu'$ be another $K$-valued
random variable, such that $\nu^*$ and $\nu'$ are i.i.d.~given $Y'$.
Then it is not difficult to argue that for any measurable function $h$,
\[
\ee\bigl\llVert \nu^*-h\bigl(Y'\bigr)\bigr\rrVert ^2
\ge\tfrac{1}{2}\ee\bigl\llVert \nu^*-\nu'\bigr\rrVert
^2.
\]
In particular, this holds for $h=g$. Thus, it suffices to show that
%
\begin{equation}
\label{s2} \ee\bigl\llVert \nu^* - \nu'\bigr\rrVert
^2 \ge C_5 t_{\mu^*}^2
\end{equation}
for some constant $C_5$.

Let $B_1$ be the ball of radius $C_6 t_{\mu^*}$ around $\nu^*$, where
$C_6$ is a constant that will be chosen later. Note that unlike $B_0$,
$B_1$ is a random set. If we can show that $\pp(\nu'\notin B_1)$ is
larger than a universal threshold, it will complete the proof of (\ref{s2}).

To prove this, the first step is to explicitly write down
\begin{eqnarray*}
\pp\bigl(\nu'\in B_1\mid Y', \nu^*\bigr)
&=& \frac{\int_{B_1} e^{-Z'\cdot(\nu^*-\nu) - {1}/{2}\llVert  \nu
^*-\nu\rrVert  ^2} \,d\rho(\nu)}{\int_{B_0} e^{-Z'\cdot(\nu^*-\nu) -
{1}/{2}\llVert  \nu^*-\nu\rrVert  ^2} \,d\rho(\nu)}.
\end{eqnarray*}
After a sequence of relatively complicated technical steps involving
concentration inequalities and second moment lower bounds, one can
produce an upper bound on the expectation of the right-hand side. The
complications arise from the fact that the right-hand side is a ratio
of random variables. When the dust settles, we get the inequality
\[
\pp\bigl(\nu'\in B_1\bigr)\le C_7 \sqrt{
\ee\bigl(\rho(B_1)\bigr)}.
\]
The proof, therefore, will be complete if we can show that $\ee(\rho
(B_1))$ can be made as small as we like by choosing $C_6$ small enough.
By the definition of $\rho$, it is clear that
\[
\rho(B_1) \le\pp\bigl(M_2\ge M_1
\mid\nu^*\bigr),
\]
where
\[
M_1 = \sup_{\mu\in B_0} Z'\cdot\bigl(\mu-
\mu^*\bigr),\qquad M_2 = \sup_{\mu\in B_1} Z'\cdot
\bigl(\mu-\mu^*\bigr).
\]
Consequently, $\ee(\rho(B_1))\le\pp(M_2\ge M_1)$. To make the
right-hand side as small as we need it to be, it makes sense to choose
$\mu^*$ such that $\ee(M_1)$ is as large as possible, and then choose
$C_6$ so small that $\ee(M_2)$ is small enough. Working out the
details of this step involves delicate technical problems. Carefully
solving these problems leads to the completion of the proof of~Theorem
\ref{admissibility}.

\section{Proofs}\label{proofs}
This section contains the proofs of all the results stated in
Sections~\ref{theory} and~\ref{applications}. We will follow a
certain notational convention about universal constants throughout this
section. Within the proof of each lemma or theorem or proposition,
$C_1, C_2,\ldots$ will denote positive universal constants. The values
of the $C_i$'s may change from one lemma to the next. On the other
hand, $c_1, c_2,\ldots$ will denote universal constants whose values
are important; once defined, they will not change.

The first goal is to prove Theorem~\ref{mainest}. We need the
following ingredient from measure concentration theory.

%
\begin{lmm}[(Cirelson, Ibragimov and Sudakov \cite{tis76})]\label{gaussmax}
Let $V_1, \ldots, V_n$ be jointly Gaussian random variables, each with
mean zero and second moment bounded above by $1$ (but not necessarily
independent). Let $M:= \max_{1\le i\le n} V_i$. Then for any $t\ge0$,
\begin{eqnarray*}
\max\bigl\{\pp\bigl(M - \ee(M) \ge t\bigr), \pp\bigl(M - \ee(M) \le-t\bigr)
\bigr\} &\le& e^{-t^2/2}.
\end{eqnarray*}
\end{lmm}
The above inequalities were proved in \cite{tis76}, although they
follow (with slightly worse constants) from the earlier papers \cite
{borell75} and \cite{sudakovtsirelson74}.

We also need a standard fact from convex geometry. A proof is included
for the sake of completeness.

%
\begin{lmm}[(Projection on to convex sets)]\label{proj}
Let $K$ be a nonempty closed convex subset of $\rr^n$. For any $x\in
\rr^n$, there is a unique point in $K$, that we call $P_K(x)$, which
is closest to $x$.
\end{lmm}

\begin{pf}
Let $s:= \inf_{y\in K}\llVert  x-y\rrVert  $. Since $K$ is nonempty, $s$ is finite.
Let $K'$ be the set of all points in $K$ that are within distance $s+1$
from $x$. This is clearly nonempty, convex and bounded. Furthermore,
since $K$ is closed, so is $K'$. The compactness of $K'$ ensures the
existence of at least one point in $K'$ that is at distance exactly $s$
from $x$. This proves the existence of a projection. Suppose now that
there are two points $y$ and $z$ in $K$ that are both at distance
exactly $s$ from $x$. Then the three points $x$, $y$ and $z$ form an
isosceles triangle with the line segment joining $y$ and $z$ as the
base. But this line segment is contained in $K$, because $K$ is convex.
Since $y\ne z$, this proves that there is a point in $K$ that is at
distance strictly less than $s$ from $x$, which is impossible.
\end{pf}

We are now ready to prove Theorem~\ref{mainest}, Corollary~\ref
{maincor} and Proposition~\ref{simple}.

\begin{pf*}{Proof of Theorem~\ref{mainest}}
Fix $\mu\in\rr^n$ and let $Y = Z+\mu$. Define two random functions
$M$ and $F$ from $[0,\infty)$ into $[-\infty, \infty)$ as
\[
M(t):= \sup_{\nu\in K\dvtx   \llVert  \nu-\mu\rrVert  \le t} Z\cdot(\nu-\mu)
\]
and
\[
F(t):= M(t) - \frac{t^2}{2},
\]
with the usual convention that the supremum of an empty set is $-\infty
$. Let $m(t):= \ee(M(t))$. Note that $\ee(F(t)) = f_\mu(t) = m(t)-t^2/2$.

Note that $M(t)$, $F(t)$, $m(t)$ and $f_\mu(t)$ are all finite if
$t\ge t_c$, and $-\infty$ if $t<t_c$. Take any $t_c\le s\le t$. Let
$\nu_1$ and $\nu_2$ be points in $K$ such that
%
\begin{equation}
\label{nucrit} \llVert \nu_1-\mu\rrVert \le s \quad\mbox{and}\quad
\llVert \nu_2-\mu\rrVert \le t.
\end{equation}
Take any $u\in[0,1]$ and let $\nu:= u\nu_1 + (1-u)\nu_2$. Then $\llVert
\nu-\mu\rrVert   \le r:= us + (1-u)t$. On the other hand,
\[
Z\cdot(\nu-\mu) = u Z\cdot(\nu_1- \mu) + (1-u) Z\cdot(
\nu_2 - \mu).
\]
Maximizing over all $\nu_1$ and $\nu_2$ satisfying (\ref{nucrit}),
this gives
%
\begin{equation}
\label{mconv} M(r)\ge u M(s)+(1-u) M(t).
\end{equation}
Thus, $M$ is a concave function of $t$. Consequently, $F$ is strictly
concave. Note that $\lim_{t\rightarrow\infty} F(t)=-\infty$, since
%
\begin{equation}
\label{mbound} M(t) \le\sup_{\nu\in\rr^n\dvtx    \llVert  \nu-\mu\rrVert  \le t} Z\cdot(\nu -\mu) = t \llVert
Z\rrVert.
\end{equation}
The strict concavity and the decay to $-\infty$ prove the existence
and uniqueness of a (random) point $t^*\in[t_c,\infty)$ where $F$ is
maximized.

Taking expectation on both sides in (\ref{mconv}) implies that $m$ is
also concave, and therefore $f_\mu$ is strictly concave. Similarly by
(\ref{mbound}), $m(t)\le t  \ee\llVert  Z\rrVert  $, which proves that $\lim_{t\rightarrow\infty} f_\mu(t)=-\infty$. Therefore, $t_\mu$ exists
and is unique.

Let $\nu^*$ be a point in $K$ that maximizes $Z\cdot(\nu-\mu)$
among all $\nu\in K$ satisfying $\llVert  \nu-\mu\rrVert  \le t^*$. Let $t_0:= \llVert
\nu^*-\mu\rrVert  $. If $t_0< t^*$, then
\[
F(t_0)\ge Z\cdot\bigl(\nu^*-\mu\bigr) - \frac{t_0^2}{2} = M
\bigl(t^*\bigr) - \frac
{t_0^2}{2} > F\bigl(t^*\bigr),
\]
which is false. Therefore, $t_0=t^*$. This shows that for any $\nu\in K$,
\begin{eqnarray*}
Z\cdot(\nu-\mu) - \frac{\llVert  \nu-\mu\rrVert  ^2}{2} &\le& F\bigl(\llVert \nu-\mu\rrVert \bigr)
\\
&\le& F\bigl(t^*\bigr) = Z\cdot\bigl(\nu^*-\mu\bigr) - \frac{\llVert  \nu^*-\mu\rrVert  ^2}{2}.
\end{eqnarray*}
Since
\[
\llVert Y-\nu\rrVert ^2 = \llVert Y-\mu\rrVert ^2 - 2
\biggl(Z\cdot(\nu-\mu) - \frac{\llVert
\nu-\mu\rrVert  ^2}{2} \biggr),
\]
this proves that $\llVert  Y-\nu\rrVert  \ge\llVert  Y-\nu^*\rrVert  $ for all $\nu\in K$.
Therefore, by the uniqueness of projection on to closed convex sets,
$\hmm= \nu^*$. In particular,
\[
\llVert \mu-\hmm\rrVert =t^*.
\]
Now note that for any $t\ge t_c$, the inequality $f_\mu(t)\le f_\mu
(t_\mu)$ may be rewritten as
%
\begin{equation}
\label{mmu1} m(t) \le m(t_\mu) + \frac{t^2-t_\mu^2}{2}.
\end{equation}
By concavity of $m$, for any $\varepsilon\in(0,1)$,
%
\begin{equation}
\label{mmu2} m\bigl((1-\varepsilon)t_\mu+ \varepsilon t\bigr) \ge(1-
\varepsilon) m(t_\mu) + \varepsilon m(t).
\end{equation}
Applying (\ref{mmu1}) to $(1-\varepsilon)t_\mu+ \varepsilon t$ instead of
$t$ gives
\[
m\bigl((1-\varepsilon)t_\mu+ \varepsilon t\bigr) \le m(t_\mu)
+ \frac{(-2\varepsilon+
\varepsilon^2)t_\mu^2 + 2(1-\varepsilon)\varepsilon t_\mu t +\varepsilon^2
t^2}{2}.
\]
Combining this inequality with (\ref{mmu2}) gives
\[
\varepsilon m(t) \le\varepsilon m(t_\mu) + \frac{(-2\varepsilon+
\varepsilon^2)t_\mu^2 + 2(1-\varepsilon)\varepsilon t_\mu t +\varepsilon^2
t^2}{2}.
\]
Dividing both sides by $\varepsilon$ and taking $\varepsilon\rightarrow0$,
we get
%
\begin{equation}
\label{mmineq} m(t) \le m(t_\mu) - t_\mu^2 +
t_\mu t,
\end{equation}
which may be rewritten as
\begin{eqnarray*}
f_\mu(t) &\le& f_\mu(t_\mu) -
\frac{(t-t_\mu)^2}{2}.
\end{eqnarray*}
Note that the above two inequalities hold even if $t<t_c$. Take any
$x>0$ and let $r_1:= t_\mu- x\sqrt{t_\mu}$ and $r_2:= t_\mu+
x\sqrt{t_\mu}$. First assume that $r_1 \ge t_c$. Then by the above inequality,
\[
\max\bigl\{f_\mu(r_1), f_\mu(r_2)
\bigr\} \le f_\mu(t_\mu) - \frac
{x^2t_\mu}{2}.
\]
%
By the concentration inequality for maxima of Gaussian random variables
(Lemma~\ref{gaussmax}), for any $t\ge0$ and $y\ge0$,
\[
\max\bigl\{\pp\bigl(F(t) \ge f_\mu(t) + y\bigr), \pp\bigl(F(t)\le
f_\mu(t)-y\bigr)\bigr\}\le e^{-y^2/2t^2}.
\]
Taking $y = x^2t_\mu/4$ and $z = f_\mu(t_\mu)-y$, a combination of
the last two displays gives the inequalities
\begin{eqnarray*}
\pp\bigl(F(r_1) \ge z\bigr) &\le&\pp\bigl(F(r_1) \ge
f_\mu(r_1) + y\bigr) \le e^{-y^2/2r_1^2},
\\
\pp\bigl(F(r_2) \ge z\bigr) &\le&\pp\bigl(F(r_2) \ge
f_\mu(r_2) + y\bigr) \le e^{-y^2/2r_2^2},
\\
\pp\bigl(F(t_\mu) \le z\bigr) &=& \pp\bigl(F(t_\mu) \le
f_\mu(t_\mu) - y\bigr) \le e^{-y^2/2t_\mu^2}.
\end{eqnarray*}
Let $E$ be the event that $F(r_1)< z$, $F(r_2)< z$ and $F(t_\mu) > z$.
By the above three inequalities,
\[
\pp\bigl(E^c\bigr) \le e^{-y^2/2r_1^2} + e^{-y^2/2r_2^2} +
e^{-y^2/2t_\mu^2}\le 3 e^{-y^2/2r_2^2}.
\]
On the other hand, by the concavity of $F$, if $E$ happens then $t^*$
must lie in the interval $(r_1,r_2)$. Together with our previous
observation that $t^*=\llVert  \mu-\hmm\rrVert  $, this completes the proof of the
theorem when $r_1\ge t_c$.

If $r_1 < t_c$, the inequality $f_\mu(r_2) \le f_\mu(t_\mu) -
x^2t_\mu/2$ is still true. Redefine $E$ to be the event that $F(r_2) <
z$ and $F(t_\mu) > z$. Then the upper bound on $\pp(E^c)$ is still
valid, and the occurrence of $E$ implies that $t^*\in[t_c,
r_2)\subseteq(r_1,r_2)$. This finishes the argument in the case $r_1 < t_c$.
\end{pf*}

\begin{pf*}{Proof of Corollary~\ref{maincor}}
Throughout this proof, $C$ denotes an arbitrary universal constant
whose value may change from line to line. First, suppose that $t_\mu
\ge1$. Then by Theorem~\ref{mainest},
\[
\pp \bigl( \bigl\llvert \llVert \hmm-\mu\rrVert - t_\mu \bigr\rrvert
\ge x\sqrt{t_\mu } \bigr) \le3 e^{-x^4/32(1+x)^2}.
\]
This shows that
\[
\ee\bigl(\llVert \hmm-\mu\rrVert -t_\mu\bigr)^2 \le
Ct_\mu,
\]
which gives the first set of inequalities. On the other hand, if $t_\mu
< 1$, then putting $z= x\sqrt{t_\mu}$, Theorem~\ref{mainest} gives
\[
\pp \bigl( \bigl\llvert \llVert \hmm-\mu\rrVert - t_\mu \bigr\rrvert
\ge z \bigr) \le3 e^{-z^4/32(t_\mu+z)^2} \le3e^{-z^4/32(1+z)^2},
\]
which gives the second inequality.
\end{pf*}

\begin{pf*}{Proof of Proposition~\ref{simple}}
The first two assertions are obvious by the strict concavity of $f_\mu
$. For the third one, observe that if $\mu\in K$, then $f_\mu(0)=0$,
and apply the second assertion.
\end{pf*}

The next goal is to prove Theorem~\ref{admissibility}. In addition to
Lemmas~\ref{gaussmax}~and~\ref{proj}, we need a few more
standard results. The first result, stated below, is called the
``Gaussian concentration inequality.''

%
\begin{lmm}[(Gaussian concentration inequality)]\label{gaussconc}
Let $Z$ be an $n$-dimen\-sional standard Gaussian random vector, and let
$f\dvtx \rr^n \rightarrow\rr$ be a function that satisfies
$\llvert  f(x)-f(y)\rrvert  \le L\llVert  x-y\rrVert  $ for all $x$ and $y$, where $L$ is a positive
constant. Then for any $\theta\in\rr$,
\[
\ee\bigl(e^{\theta(f(Z)-\ee(f(Z)))}\bigr) \le e^{L^2\theta^2/2}.
\]
Consequently, for any $t\ge0$,
\[
\max\bigl\{\pp\bigl(f(Z)-\ee\bigl(f(Z)\bigr) \ge t\bigr), \pp\bigl(f(Z)-\ee
\bigl(f(Z)\bigr)\le-t\bigr)\bigr\} \le e^{-t^2/2L^2}.
\]
\end{lmm}
This famous result possibly appeared for the first time as an implied
consequence of the theorems in \cite
{borell75,sudakovtsirelson74,tis76}. For a simple proof, originally
appearing in~\cite{tis76}, see the argument following equation (2.35)
in \cite{ledoux01}.

We also need the fact that the projection $P_K$ on to a closed convex
set is a contraction with respect to the Euclidean norm. This, again,
is quite standard but we provide a short proof for the sake of completeness.

%
\begin{lmm}\label{contrac}
For any closed convex set $K\subseteq\rr^n$ and any $x,y$, $\llVert  P_K(x)-P_K(y)\rrVert  \le\llVert  x-y\rrVert  $.
\end{lmm}

\begin{pf}
Let $z=P_K(x)$ and $w = P_K(y)$. If $z=w$, there is nothing to prove.
So assume that $z\ne w$. Let $S$ be the line segment joining $z$ and
$w$. Then $S$ is entirely contained in $K$. Let $H_1$ be the hyperplane
passing through $z$ that is orthogonal to~$S$, and $H_2$ be the
hyperplane passing through $w$ that is orthogonal to~$S$.

The hyperplane $H_1$ divides $\rr^n\setminus H_1$ into two open
half-spaces, one of which contains $w$. If $x$ belongs to the
half-space that contains $w$, then there is a point on $S$ that is
closer to $x$ than $z$. This is impossible. Similarly, $H_2$ divides
$\rr^n\setminus H_2$ into two open half-spaces, and $y$ cannot belong
to the one that contains $z$. Therefore, both of the parallel
hyperplanes $H_1$ and $H_2$ must lie between $x$ and $y$. This proves
that $\llVert  x-y\rrVert  \ge$ the distance between $H_1$ and $H_2$, which is equal
to $\llVert  z-w\rrVert  $.
\end{pf}
Finally, we need the so-called ``second moment inequality,'' also known
as the ``Paley--Zygmund inequality.''

%
\begin{lmm}[(Second moment inequality)]\label{pzlmm}
If $X$ is a nonnegative random variable with $\ee(X)>0$ and finite
second moment, then for any $a\in[0, \ee(X)]$,
\[
\pp(X> a) \ge\frac{(\ee(X)-a)^2}{\ee(X^2)}.
\]
\end{lmm}
The proof of this standard inequality may be found in graduate
probability text books such as \cite{durrett}.

We now embark on the proof of Theorem~\ref{admissibility}. Several
preparatory lemmas are required.

\begin{lmm}\label{interval}
Let $K\subseteq\rr^n$ be a line segment of length $l$. Let $Z$ be an
$n$-dimensional standard Gaussian random vector. Let $f\dvtx \rr^n
\rightarrow\rr^n$ be any Borel measurable map. Then there exists $\mu
\in K$ such that
\[
\ee\bigl\llVert f(Z+\mu) - \mu\bigr\rrVert ^2 \ge c_1
\min\bigl\{l^2, 4\bigr\},
\]
where $c_1$ is a positive universal constant.
\end{lmm}

\begin{pf}
Let $\lambda$ denote the uniform distribution on $K$. Let $\nu$ be
point chosen uniformly at random from $K$. Let $Y:= Z+\nu$. Given $Y$
and $\nu$, let $\nu'$ be drawn from the posterior distribution of
$\nu$ given $Y$. Explicitly, if $\theta$ denotes the joint law of
$(\nu, Y, \nu')$, then
%
\begin{equation}
\label{ynu} d\theta\bigl(\mu, y, \mu'\bigr) = \frac{e^{-1/2\llVert  y-\mu'\rrVert  ^2
}e^{-1/2\llVert  y-\mu\rrVert  ^2 }}{(2\pi)^{n/2}\int_{K} e^{-1/2\llVert  y-x\rrVert  ^2} \,d\lambda(x)}
\,d\lambda(\mu) \,dy \,d\lambda\bigl(\mu '\bigr).
\end{equation}
The above expression clearly shows that $\nu$ and $\nu'$ are i.i.d. given $Y$ and, therefore,
\begin{eqnarray}\label{fynu}
\ee \bigl(\bigl\llVert f(Y) - \nu\bigr\rrVert ^2 \mid Y\bigr) &\ge&
\ee \bigl(\bigl\llVert \ee(\nu\mid Y) - \nu\bigr\rrVert ^2 \mid Y
\bigr)
\nonumber\\[-8pt]\\[-8pt]\nonumber
&=& \tfrac{1}{2} \ee \bigl(\bigl\llVert \nu' - \nu\bigr\rrVert
^2 \mid Y\bigr).
\nonumber
\end{eqnarray}
[In the above display, $\ee(\nu\mid Y)$ denotes the random vector
whose $i$th coordinate is $\ee(\nu_i\mid Y)$. The inequality in the
first line is simply a consequence of the fact that for any random
variable $X$, $\ee(X-a)^2$ is minimized when $a = \ee(X)$.]
Next, let
\[
M:= \sup_{x\in K} \bigl\llvert Z\cdot(x-\nu)\bigr\rrvert.
\]
Note that since the function being maximized is convex and the set $K$
is a line segment, therefore the maximum is necessarily attained at one
of the endpoints of the line segment $K$. From this, it easy to see that
%
\begin{equation}
\label{e2m} \ee\bigl(e^{2M}\bigr) \le C_1e^{C_1l^2}.
\end{equation}
Take any $\varepsilon>0$. Given $Y$ and $\nu$, let $I$ denote the set of
all points in $K$ that are within distance $\varepsilon$ from $\nu$.
Then by (\ref{ynu}),
\begin{eqnarray*}
\pp\bigl(\bigl\llVert \nu' - \nu\bigr\rrVert \le\varepsilon\mid Y,
\nu\bigr) &=& \frac{\int_I
e^{Z\cdot(x-\nu) - {1}/{2}\llVert  x-\nu\rrVert  ^2} \,d\lambda(x)}{\int_K
e^{Z\cdot(x-\nu) - {1}/{2}\llVert  x-\nu\rrVert  ^2} \,d\lambda(x)}
\le\frac{2e^{2M + l^2} \varepsilon}{l}.
\end{eqnarray*}
Taking expectation and applying (\ref{e2m}), we get
\[
\pp\bigl(\bigl\llVert \nu'-\nu\bigr\rrVert \le\varepsilon\bigr) \le
\frac{C_2e^{C_2l^2}\varepsilon
}{l},
\]
and, therefore,
\[
\ee\bigl\llVert \nu'-\nu\bigr\rrVert ^2 \ge
\varepsilon^2 \pp\bigl(\bigl\llVert \nu'-\nu\bigr\rrVert >
\varepsilon \bigr)\ge\varepsilon^2 \biggl(1-\frac{C_2e^{C_2l^2}\varepsilon}{l} \biggr).
\]
If $l\le2$, then combined with (\ref{fynu}) and taking $\varepsilon=
C_3 l$ for some small enough $C_3$, this proves that
\[
\ee\bigl\llVert f(Z+\nu)-\nu\bigr\rrVert ^2 \ge C_4
l^2.
\]
In particular, there exists $\mu\in K$ such that
\[
\ee\bigl\llVert f(Z+\mu)-\mu\bigr\rrVert ^2 \ge C_4
l^2.
\]
If $l> 2$, then choose a subinterval $K'\subseteq K$ of length $\le2$
and work with $K'$ instead of $K$.
\end{pf}

%
\begin{lmm}\label{tlowlmm}
There is a positive universal constant $c_2$ such that following is
true. Let $K$ be a closed convex subset of $\rr^n$ with diameter $\ge
2$. For each $\mu\in K$, let $t_\mu$ be defined as in Theorem~\ref
{mainest}. Then $t_\mu\ge c_2 n^{-1/2}$ for all $\mu\in K$.
\end{lmm}

\begin{pf}
Take any $\mu\in K$. Since the diameter of $K$ is $\ge$2, there
exists $\nu\in K$ such that $\llVert  \nu-\mu\rrVert  \ge1$. By the convexity of
$K$, this implies that there exists $\nu\in K$ such that $\llVert  \nu-\mu\rrVert
=1$. For each $t\in[0,1]$ let $\nu_t:= (1-t)\mu+ t\nu$. Then $\nu
_t\in K$ and $\llVert  \nu_t-\mu\rrVert  =t$. Therefore, there exists positive
$C_1$ and $C_2$ such that if $t\le C_1$ then
\begin{eqnarray*}
f_\mu(t) &\ge&\ee \bigl(\max\bigl\{0, Z\cdot(\nu_t-\mu)
\bigr\} \bigr)-\frac
{t^2}{2} \ge C_2t.
\end{eqnarray*}
On the other hand, by (\ref{mbound}),
\[
f_\mu(t) \le C_3t\sqrt{n}.
\]
Thus, with $C_4:= C_1C_2/C_3$,
\[
f_\mu\bigl(C_4n^{-1/2}\bigr) \le
C_2C_1 \le f_\mu(C_1).
\]
Taking $C_3$ large enough, we have $C_4n^{-1/2}< C_1$. By Proposition
\ref{simple}, this shows that $t_\mu\ge C_4 n^{-1/2}$.
\end{pf}

\begin{lmm}\label{tmubds}
Let $K$ be a closed convex subset of $\rr^n$ and let $t_\mu$ be
defined as in Theorem~\ref{mainest}. Then for any $\mu,\nu\in K$
such that $\llVert  \mu-\nu\rrVert  \le t_\mu/24$,
\[
\frac{11t_\mu}{24}\le t_\nu\le\frac{50t_\mu}{24}.
\]
\end{lmm}

\begin{pf}
If $t_\mu=0$ there is nothing to prove. So assume that $t_\mu> 0$.
For any $\gamma\in K$ and $t\ge0$, let
%
\begin{equation}
\label{bmudef} B(\gamma, t):= \bigl\{\gamma'\in K\dvtx  \bigl\llVert
\gamma'-\gamma\bigr\rrVert \le t\bigr\}
\end{equation}
and
%
\begin{equation}
\label{mudef} m_\gamma(t):= \ee \Bigl(\sup_{\gamma'\in B(\gamma, t)} Z
\cdot \bigl(\gamma'-\gamma\bigr) \Bigr) = \ee \Bigl(\sup
_{\gamma'\in B(\gamma,
t)} Z\cdot\gamma' \Bigr).
\end{equation}
Let $B_0:= B(\mu, r)$, where $r:= t_{\mu}/24$. Take any $\nu\in
B_0$. Note that for any positive integer $k$,
\[
B\bigl(\mu, (k-1)r\bigr)\subseteq B(\nu, kr) \subseteq B\bigl(\mu, (k+1)r\bigr),
\]
and, therefore,
%
\begin{equation}
\label{muk22} m_{\mu}\bigl((k-1)r\bigr)\le m_{\nu}(kr)\le
m_{\mu}\bigl((k+1)r\bigr).
\end{equation}
Applying (\ref{muk22}) with $k=11$ gives
\[
m_\nu(11r) \le m_{\mu}(12r) = m_\mu(t_\mu/2),
\]
and with $k=25$, we get
\[
m_\nu(25 r) \ge m_{\mu} (24r)= m_{\mu}(t_{\mu}).
\]
Therefore, by the inequality (\ref{mmineq}) from the proof of Theorem
\ref{mainest},
\begin{eqnarray*}
f_\nu(25r) - f_\nu(11 r) &=& m_\nu(25r) -
m_\nu(11 r) - \frac
{(25^2-11^2) r^2}{2}
\\
&\ge& m_\mu(t_\mu)-m_\mu(t_\mu/2)-
\frac{252t_\mu^2}{576}\ge\frac
{t_\mu^2}{2} - \frac{7t_\mu^2}{16}\ge0.
\end{eqnarray*}
Therefore, by Proposition~\ref{simple},
\begin{eqnarray*}
t_\nu&\ge&11 r = \frac{11t_\mu}{24}.
\end{eqnarray*}
Next, note that by (\ref{muk22}),
\begin{eqnarray*}
f_\nu(50 r ) - f_\nu(25 r) &=& m_\nu(50 r)-
m_\nu(25 r) -\frac{1875
r^2}{2}
\\
&\le& m_{\mu}(51 r) - m_{\mu}(24 r) - \frac{1875 r^2}{2}.
\end{eqnarray*}
By the inequality (\ref{mmineq}),
\begin{eqnarray*}
m_{\mu}(51 r) - m_{\mu}(24 r) &=& m_{\mu}(27 r +
t_{\mu}) - m_{\mu
}(t_{\mu}) \le27 r
t_{\mu} = 648 r^2.
\end{eqnarray*}
Combining the last two displays gives
\begin{eqnarray*}
f_\nu(50 r ) - f_\nu(25 r) &\le&648 r^2 -
\frac{1875r^2}{2}\le0.
\end{eqnarray*}
By Proposition~\ref{simple}, this proves that $t_\nu\le50r$.
\end{pf}

We are now ready to prove Theorem~\ref{admissibility}.

\begin{pf*}{Proof of Theorem~\ref{admissibility}}
First, suppose that $l:=\operatorname{diam}(K) \le2$. Choose a~line
segment $I\subseteq K$ of length $l$. By Lemma~\ref{interval}, there
exists $\mu_0\in I$ such that $\ee\llVert  g(Z+\mu_0)-\mu_0\rrVert  ^2 \ge c_1
l^2$. But $\ee\llVert  P_K(Z+\mu_0)-\mu_0\rrVert  ^2 \le l^2$, since any two
elements of $K$ are within distance $l$ of each other. This completes
the proof of the theorem when $\operatorname{diam}(K)\le2$. For the
rest of the proof, assume that
$\operatorname{diam}(K) > 2$.

For $\mu\in K$ and $t\ge0$, let $m_\mu(t)$ be defined as in (\ref
{mudef}). Since $\operatorname{diam}(K)>2$, Lemma~\ref{tlowlmm}
implies that for all $\mu\in K$,
%
\begin{equation}
\label{tlow} t_\mu\ge c_2 n^{-1/2}.
\end{equation}
Let
\[
s:= \sup_{\mu\in K}m_\mu\bigl( 10^{-3}t_\mu
\bigr).
\]
Then there exists at least one point $\mu^*\in K$ such that
\[
m_{\mu^*}\bigl(10^{-3}t_{\mu^*}\bigr) \ge s -
\frac{c_2^2}{10^6 n}.
\]
For $\nu\in K$ and $t\ge0$ let $B(\nu, t)$ be defined as in (\ref
{bmudef}). Let $B_0:= B(\mu^*, r)$, where $r:= 10^{-3}t_{\mu^*}$.
Lemma~\ref{tmubds} implies that for all $\nu\in B_0$,
%
\begin{equation}
\label{tlow23} \frac{11t_{\mu^*}}{24}\le t_\nu\le\frac{50t_{\mu^*}}{24}.
\end{equation}
Define a probability measure $\rho$ on $B_0$ as follows. Let $\nu^*$
be the point that maximizes $Z\cdot(\nu-\mu^*)$ among all $\nu\in
B_0$. If there are more than one such points, take the one that is the
least in the lexicographic ordering (it is easy to prove that there is
a least element since the set of maximizers is closed). Let $\rho$ be
the law of $\nu^*$. Let $Z'$ be a standard Gaussian random vector,
independent of $Z$. Let $Y':= Z' + \nu^*$. Let $\rho'$ be the
conditional distribution of $\nu^*$ given $Y'$. It is easy to see that
\[
d\rho'(\nu) = L^{-1} e^{-{1}/{2}\llVert  Y'-\nu\rrVert  ^2} \,d\rho(\nu),\qquad \nu\in
B_0,
\]
where
\[
L:= \int_{B_0} e^{-{1}/{2}\llVert  Y'-\nu\rrVert  ^2} \,d\rho(\nu).
\]
Given $Y'$ and $\nu^*$, let $\nu'$ be a random point generated from
the distribution $\rho'$. Then $\nu'$ and $\nu^*$ are conditionally
i.i.d. given $Y'$, as is evident from the joint law $\theta$ of the
triple $(\nu^*,Y', \nu')$:
%
\begin{eqnarray}
d\theta(\nu_1,y, \nu_2) &=& \frac{e^{-1/2\llVert  y-\nu_2\rrVert
^2}e^{-1/2\llVert  y-\nu_1\rrVert  ^2}}{(2\pi)^{n/2}\int_{B_0} e^{-1/2\llVert  y-\nu\rrVert  ^2} \,d\rho(\nu)} \,d\rho(
\nu_1) \,dy \,d\rho(\nu _2),\nonumber
\\
\eqntext{\mbox{where } \nu_1,\nu_2\in B_0
\mbox{ and } y\in\rr^n.}
\end{eqnarray}
Let $\ee(\nu^*\mid Y')$ be the random vector whose $i$th coordinate
is $\ee(\nu^*_i\mid Y')$ and $g$ be an arbitrary Borel measurable map
from $\rr^n$ into itself, as in the statement of Theorem~\ref
{admissibility}. Then
\begin{eqnarray*}
\ee\bigl(\bigl\llVert \nu^*-\nu'\bigr\rrVert ^2 \mid
Y'\bigr) &=& 2 \ee\bigl(\bigl\llVert \nu^* - \ee\bigl(\nu ^*\mid
Y'\bigr)\bigr\rrVert ^2 \mid Y'\bigr)
\\
&\le&2 \ee\bigl(\bigl\llVert \nu^* - g\bigl(Y'\bigr)\bigr\rrVert
^2 \mid Y'\bigr).
\end{eqnarray*}
Thus,
%
\begin{eqnarray}
\label{nnu} \ee\bigl\llVert \nu^* - \nu'\bigr\rrVert
^2 &\le&2 \ee\bigl\llVert \nu^* - g\bigl(Y'\bigr)\bigr
\rrVert ^2.
\end{eqnarray}
Let $B_1$ denote the (random) set $B(\nu^*, 10^{-3}r)\cap B_0$. Then
%
\begin{eqnarray}\label{nub1}
\pp\bigl(\nu'\in B_1\mid Y', \nu^*\bigr)
&=& L^{-1}\int_{B_1} e^{-1/2\llVert  Y'-\nu\rrVert  ^2} \,d\rho(\nu)
\nonumber\\[-8pt]\\[-8pt]\nonumber
&=& \frac{\int_{B_1} e^{-Z'\cdot(\nu^*-\nu) - 1/2\llVert  \nu
^*-\nu\rrVert  ^2} \,d\rho(\nu)}{\int_{B_0} e^{-Z'\cdot(\nu^*-\nu) -
1/2\llVert  \nu^*-\nu\rrVert  ^2} \,d\rho(\nu)}.
\end{eqnarray}
Let $L_1$ and $L_2$ denote the numerator and the denominator in the
last expression. First, note that
%
\begin{eqnarray}\label{nub2}
\ee\bigl(L_1^2\mid\nu^*\bigr) &\le&\int
_{B_1}\ee\bigl(e^{-2Z'\cdot(\nu^*-\nu)
- \llVert  \nu^*-\nu\rrVert  ^2}\mid\nu^*\bigr) \,d\rho(
\nu)
\nonumber\\[-8pt]\\[-8pt]\nonumber
&=& \int_{B_1}e^{\llVert  \nu^*-\nu\rrVert  ^2} \,d\rho(\nu)\le e^{10^{-6}r^2}
\rho(B_1).
\end{eqnarray}
Next, note that $\ee(L_2\mid\nu^*) = 1$, and
\begin{eqnarray*}
&& \ee\bigl(L_2^2 \mid\nu^*\bigr)
\\
&&\qquad = \int_{B_0}\int_{B_0} \ee
\bigl(e^{-Z'\cdot((\nu^*-\nu_1) + (\nu
^*-\nu_2)) - 1/2(\llVert  \nu^*-\nu_1\rrVert  ^2+\llVert  \nu^*-\nu_2\rrVert  ^2)}\bigr) \,d\rho(\nu_1) \,d\rho(\nu_2)
\\
&&\qquad = \int_{B_0}\int_{B_0} e^{(\nu^*-\nu_1) \cdot(\nu^*-\nu_2)}
\,d\rho(\nu_1) \,d\rho(\nu_2)\le e^{4r^2}.
\end{eqnarray*}
Therefore by the second moment inequality (Lemma~\ref{pzlmm}),
%
\begin{equation}
\label{l21} \pp\bigl(L_2> 1/2\mid\nu^*\bigr) \ge\frac{(\ee(L_2\mid\nu^*))^2}{4  \ee
(L_2^2\mid\nu^*)}
\ge\frac{1}{4}e^{-4r^2}.
\end{equation}
Now note that, by a slight abuse of notation,
\begin{eqnarray*}
\frac{\partial}{\partial Z_i'} \log L_2 &=& -\frac{1}{L_2}\int
_{B_0}\bigl(\nu^*_i-\nu_i\bigr)
e^{-Z'\cdot(\nu^*-\nu) - 1/2\llVert  \nu
^*-\nu\rrVert  ^2} \,d\rho(\nu).
\end{eqnarray*}
Consequently,
%
\begin{eqnarray}\label{gradest}
\sum_{i=1}^n \biggl(\frac{\partial}{\partial Z_i'}
\log L_2 \biggr)^2 &\le&\frac{\int_{B_0}\llVert  \nu^*-\nu\rrVert  ^2  e^{-Z'\cdot(\nu^*-\nu
) - 1/2\llVert  \nu^*-\nu\rrVert  ^2} \,d\rho(\nu)}{ \int_{B_0}
e^{-Z'\cdot(\nu^*-\nu) - 1/2\llVert  \nu^*-\nu\rrVert  ^2} \,d\rho(\nu)}
\nonumber\\[-8pt]\\[-8pt]\nonumber
&\le&4r^2.
\end{eqnarray}
Therefore, by the Gaussian concentration inequality (Lemma~\ref{gaussconc}), for any $x\ge0$,
%
\begin{equation}
\label{concapp} \pp\bigl(\log L_2 \ge\ee\bigl(\log L_2\mid
\nu^*\bigr) + x\mid\nu^*\bigr) \le e^{-x^2/8r^2}.
\end{equation}
Now suppose that $4r^2 > \log4$, or in other words,
%
\begin{equation}
\label{problem} t_{\mu^*} > 500\sqrt{2\log2}.
\end{equation}
Under the above condition, taking $x = 8r^2$ in (\ref{concapp}) gives
\[
\pp\bigl(\log L_2 \ge\ee\bigl(\log L_2\mid\nu^*\bigr) +
8r^2\mid\nu^*\bigr) \le e^{-8r^2}< \tfrac{1}{4}e^{-4r^2}.
\]
Comparing this with (\ref{l21}), we realize that under (\ref
{problem}), it must be true that
\[
\ee\bigl(\log L_2\mid\nu^*\bigr) \ge-8r^2 - \log2.
\]
Therefore, if (\ref{problem}) holds, then
\begin{eqnarray*}
\ee\bigl(L_2^{-2}\mid\nu^*\bigr) &=& e^{-2  \ee(\log L_2\mid\nu^*)} \ee
\bigl(e^{-2(\log L_2 - \ee(\log L_2\mid\nu^*))}\mid\nu^*\bigr)
\\
&\le&4 e^{16r^2} \ee\bigl(e^{-2(\log L_2 - \ee(\log L_2\mid\nu^*))}\mid \nu^*\bigr).
\end{eqnarray*}
But by the Gaussian concentration inequality (Lemma~\ref{gaussconc})
and the estimate~(\ref{gradest}),
\[
\ee\bigl(e^{-2(\log L_2 - \ee(\log L_2\mid\nu^*))}\mid\nu^*\bigr)\le e^{8r^2}.
\]
Combining the last two displays gives
%
\begin{eqnarray}
\label{nub3} \ee\bigl(L_2^{-2}\mid\nu^*\bigr) &\le& 4e^{24r^2}.
\end{eqnarray}
By (\ref{nub1}), (\ref{nub2}) and (\ref{nub3}), we see that under
condition (\ref{problem}),
%
\begin{eqnarray}\label{mm1}
\pp\bigl(\nu'\in B_1\mid\nu^*\bigr) &=& \ee
\bigl(L_1L_2^{-1}\mid\nu^*\bigr)\nonumber
\\
&\le& \bigl(\ee\bigl(L_1^2\mid\nu^*\bigr) \ee
\bigl(L_2^{-2}\mid\nu^*\bigr) \bigr)^{1/2}
\\
&\le&2 e^{13r^2}\sqrt{\rho(B_1)}.\nonumber
\end{eqnarray}
Define
\begin{eqnarray*}
M_1&:=& \sup_{\nu\in B_0} Z'\cdot\bigl(\nu-
\mu^*\bigr),\qquad
M_2:= \sup_{\nu\in B_1} Z'\cdot
\bigl(\nu-\nu^*\bigr),
\\
M_3&:=& Z'\cdot\bigl(\nu^*-\mu ^*
\bigr).
\end{eqnarray*}
The basic fact, easy to see, is that
%
\begin{eqnarray}\label{mm2}
\rho(B_1) &\le&\pp \Bigl(\sup_{\nu\in B_1}
Z'\cdot\bigl(\nu-\mu^*\bigr) \ge\sup_{\nu\in B_0}
Z'\cdot\bigl(\nu-\mu^*\bigr)\mid\nu^* \Bigr)
\nonumber\\[-8pt]\\[-8pt]\nonumber
&\le&\pp\bigl(M_2+M_3 \ge M_1 \mid\nu^*
\bigr).
\end{eqnarray}
Having understood this, note that by the definitions of $\mu^*$ and
$s$ and the lower bounds (\ref{tlow}) and~(\ref{tlow23}),
%
\begin{eqnarray}\label{m11}
\ee\bigl(M_1\mid\nu^*\bigr) &=& m_{\mu^*}
\bigl(10^{-3}t_{\mu^*}\bigr) \ge s - \frac
{c_2^2}{10^6n}
\nonumber\\[-8pt]\\[-8pt]\nonumber
&\ge& m_{\nu^*}\bigl(10^{-3}t_{\nu^*}\bigr) -
\frac{t_{\mu^*}^2}{10^6} \ge m_{\nu^*}(11r/24) - r^2.
\end{eqnarray}
On the other hand,
%
\begin{eqnarray}
\ee\bigl(M_2\mid\nu^*\bigr) &\le& m_{\nu^*}
\bigl(10^{-3}r\bigr). \label{m12}
\end{eqnarray}
Let $\delta:= 11r/24 - 10^{-3} r$. By the concavity of $m_{\nu^*}$,
and the inequalities (\ref{mmineq}) and~(\ref{tlow23}),
%
\begin{eqnarray}\label{m13}
&& m_{\nu^*}(11r/24) - m_{\nu^*}\bigl(10^{-3} r\bigr)\nonumber
\\
&&\qquad  =
m_{\nu^*}(11r/24) - m_{\nu^*}(11r/24 - \delta)
\nonumber\\[-8pt]\\[-8pt]\nonumber
&&\qquad \ge m_{\nu^*}(t_{\nu^*}) - m_{\nu^*}(t_{\nu^*}
- \delta)
\\
&&\qquad \ge t_{\nu^*}\delta\ge\frac{11t_{\mu^*}\delta}{24}\ge\frac
{110t_{\mu^*}r}{24^2} \ge100
r^2.\nonumber
\end{eqnarray}
By (\ref{m11}), (\ref{m12}) and (\ref{m13}), we see that
\[
\ee\bigl(M_1\mid\nu^*\bigr) - \ee\bigl(M_2\mid\nu^*
\bigr) \ge99r^2.
\]
Let $x = 33r^2$. Then by the above inequality,
\begin{eqnarray*}
&&\pp\bigl(M_2+M_3\ge M_1 \mid\nu^*\bigr)
\\
&&\qquad \le\pp\bigl(M_1\le\ee\bigl(M_1\mid\nu^*\bigr)-x\mid
\nu^*\bigr)
\\
&&\quad\qquad{} + \pp\bigl(M_2\ge\ee\bigl(M_1\mid\nu^*\bigr)-2x
\mid\nu^*\bigr) + \pp\bigl(M_3\ge x\mid\nu^*\bigr)
\\
&&\qquad \le\pp\bigl(M_1\le\ee\bigl(M_1\mid\nu^*\bigr)-x\mid
\nu^*\bigr)
\\
&&\quad\qquad{} + \pp\bigl(M_2\ge\ee\bigl(M_2\mid\nu^*\bigr)+ x
\mid\nu^*\bigr) + \pp\bigl(M_3\ge x\mid\nu^*\bigr).
\end{eqnarray*}
By the concentration inequality for Gaussian maxima (Lemma~\ref{gaussmax}) and the fact that $\ee(M_3\mid\nu^*)=0$, this shows that
\begin{eqnarray*}
\pp\bigl(M_2+M_3\ge M_1 \mid\nu^*\bigr) &\le& e^{-x^2/2r^2} + e^{-x^2
/2(10^{-3}r)^2} + e^{-x^2/2r^2}
\\
&\le&3 \exp\bigl(-500r^2\bigr).
\end{eqnarray*}
Combined with (\ref{mm1}) and (\ref{mm2}), this shows that if (\ref
{problem}) holds, then
\[
\pp\bigl(\nu'\in B_1\mid\nu^*\bigr) \le
C_1 \exp\bigl(-C_2 t_{\mu^*}^2\bigr).
\]
Therefore, there is a universal constant $C_3\ge500\sqrt{2\log2}$
such that if $t_{\mu^*}\ge C_3$, then
\begin{eqnarray*}
\ee\bigl\llVert \nu' - \nu^*\bigr\rrVert ^2 &\ge&
\bigl(10^{-3}r\bigr)^2 \pp\bigl(\nu'\notin
B_1\bigr) \ge C_4 t_{\mu^*}^2,
\end{eqnarray*}
and so by (\ref{nnu}),
\[
\ee\bigl\llVert \nu^*- g\bigl(Z'+\nu^*\bigr)\bigr\rrVert
^2 \ge C_5 t_{\mu^*}^2.
\]
Since
\[
\ee\bigl\llVert \nu^*- g\bigl(Z'+\nu^*\bigr)\bigr\rrVert
^2 = \int_{B_0} \ee\bigl\llVert \mu- g
\bigl(Z'+\mu\bigr)\bigr\rrVert ^2 \,d\rho(\mu),
\]
this shows that there exists $\mu_0\in B_0$ such that
\[
\ee\bigl\llVert \mu_0- g(Z+\mu_0)\bigr\rrVert
^2 \ge C_5t_{\mu^*}^2.
\]
By (\ref{tlow23}), $t_{\mu^*}\ge24t_{\mu_0}/50$. On the other hand
if $t_{\mu^*}\ge C_3$, then by (\ref{tlow23}), $t_{\mu_0} \ge
11t_{\mu^*}/24 \ge200\sqrt{2\log2}$. Therefore, by Corollary~\ref{maincor},
\[
\ee\bigl\llVert \mu_0- g(Z+\mu_0)\bigr\rrVert
^2 \ge C_6t_{\mu_0}^2\ge
C_7\ee\bigl\llVert \mu_0- P_K(Z+
\mu_0)\bigr\rrVert ^2.
\]
This completes the proof of the theorem when $t_{\mu^*}\ge C_3$ and
$\operatorname{diam}(K)> 2$.

Suppose now that $t_{\mu^*}< C_3$ and $\operatorname{diam}(K)> 2$.
For each $\mu$, let
\[
l_\mu^2:= \ee\bigl\llVert P_K(Z+\mu)-\mu
\bigr\rrVert ^2.
\]
Then by Corollary~\ref{maincor}, $l_{\mu^*}\le C_8$. Let $I$ be a
line segment in $K$ of length $1$, with one endpoint at $\mu^*$. By
Lemma~\ref{interval}, there exists $\mu_0\in I$ such that
%
\begin{equation}
\label{lasteq} \ee\bigl\llVert g(Z+\mu_0)-\mu_0\bigr
\rrVert ^2 \ge c_1.
\end{equation}
On the other hand, by Lemma~\ref{contrac},
\begin{eqnarray*}
&& \bigl\llVert P_K(Z+\mu_0)-\mu_0\bigr
\rrVert
\\
&&\qquad  \le\bigl\llVert P_K(Z+\mu_0)-
P_K\bigl(Z+\mu^*\bigr)\bigr\rrVert
+ \bigl\llVert P_K\bigl(Z+\mu^*\bigr) - \mu^*\bigr\rrVert +
\bigl\llVert \mu^*-\mu_0\bigr\rrVert
\\
&&\qquad \le\bigl\llVert P_K\bigl(Z+\mu^*\bigr) - \mu^*\bigr\rrVert + 2
\bigl\llVert \mu^*-\mu_0\bigr\rrVert
\\
&&\qquad \le\bigl\llVert P_K\bigl(Z+\mu^*\bigr) - \mu^*\bigr\rrVert + 2.
\end{eqnarray*}
Consequently,
\begin{eqnarray*}
\ee\bigl\llVert P_K(Z+\mu_0)-\mu_0
\bigr\rrVert ^2 &\le&2l_{\mu^*}^2 + 8\le
C_9.
\end{eqnarray*}
Together with (\ref{lasteq}), this completes the proof of the theorem
when $t_{\mu^*}< C_3$ and $\operatorname{diam}(K)>2$.
\end{pf*}

The next goal is to prove Proposition~\ref{counter1}. The proof is a
simple consequence of Proposition~\ref{simple}. We just have to carry
out some computations to verify the conditions of Proposition~\ref{simple}.

\begin{pf*}{Proof of Proposition~\ref{counter1}}
We have to first prove that the set $K$ is closed and convex. It is
obviously closed, and it is convex because for any $\alpha, \alpha
'\in[0,1]$ and $\theta_i, \theta_i'\in[-1,1]$,
\begin{eqnarray*}
&& t \bigl(\alpha n^{-1/4} + \alpha\theta_i n^{-1/2}
\bigr) + (1-t) \bigl(\alpha' n^{-1/4} +
\alpha' \theta'_i n^{-1/2}\bigr)
\\
&&\qquad = \alpha_t n^{-1/4} + \alpha_t
\theta_{i,t} n^{-1/2},
\end{eqnarray*}
where
\[
\alpha_t = t\alpha+ (1-t)\alpha'\in[0,1]
\]
and
\[
\theta_{i,t} = \frac{t\alpha\theta_i + (1-t)\alpha'\theta
'_i}{t\alpha+ (1-t)\alpha'}\in[-1,1].
\]
Let $\overline{Y}:= \sum_{i=1}^n Y_i/n$, so that the components of
$\tm$ are all equal to $\overline{Y}$. Defining $\bar{\mu} = \sum_{i=1}^n \mu_i/n$ and $\bar{\theta} = \sum_{i=1}^n \theta_i/n$, we have
\begin{eqnarray*}
\ee(\tm_i-\mu_i)^2 &=& \var(
\tm_i) + (\bar{\mu}-\mu_i)^2
= \frac{1+ \alpha^2 (\theta_i-\bar{\theta})^2}{n}\le\frac
{5}{n}.
\end{eqnarray*}
Therefore,
\[
\ee\llVert \tm-\mu\rrVert ^2 \le5,
\]
which proves one part of the proposition.

Next, let $\mu= (0,0,\ldots, 0)$. Take any $t\ge0$ and any $\nu\in
K$ such that $\llVert  \nu-\mu\rrVert   \le t$. Suppose that
\[
\nu_i:= \alpha n^{-1/4} + \alpha\theta_i
n^{-1/2},
\]
where $\alpha\in[0,1]$ and $\theta_i\in[-1,1]$. Note that
\begin{eqnarray*}
\llVert \nu-\mu\rrVert ^2 &\ge&\alpha^2\sqrt{n}
\bigl(1-n^{-1/4}\bigr)^2.
\end{eqnarray*}
Therefore, $\alpha\le tn^{-1/4}/(1-n^{-1/4})$. Since
\begin{eqnarray*}
Z\cdot(\nu-\mu) &=& \alpha n^{-1/4}\sum_{i=1}^n
Z_i + \alpha n^{-1/2}\sum_{i=1}^n
Z_i \theta_i,
\end{eqnarray*}
this proves that
\[
\sup_{\nu\in K\dvtx   \llVert  \nu-\mu\rrVert  \le t} Z\cdot(\nu-\mu) \le \frac{tn^{-1/2}}{1-n^{-1/4}}\Biggl
\llvert \sum_{i=1}^nZ_i
\Biggr\rrvert + \frac
{tn^{-3/4}}{1-n^{-1/4}}\sum_{i=1}^n
\llvert Z_i\rrvert.
\]
Consequently,
%
\begin{equation}
\label{fmu1} f_\mu(t)\le C_1tn^{1/4} -
\frac{t^2}{2}.
\end{equation}
On the other hand, if $t\le n^{1/4}$, then taking $\theta
_i=\operatorname{sign}(Z_i)$ and $\alpha= tn^{-1/4}/2$, we get $\llVert  \nu
-\mu\rrVert  \le t$ and
\[
Z\cdot(\nu-\mu) = \frac{tn^{-1/2}}{2}\sum_{i=1}^nZ_i
+ \frac
{tn^{-3/4}}{2}\sum_{i=1}^n \llvert
Z_i\rrvert,
\]
proving that
%
\begin{equation}
\label{fmu2} f_\mu(t)\ge C_2tn^{1/4} -
\frac{t^2}{2}.
\end{equation}
Without loss of generality, assume that $C_2< 1 < C_1$. Let
\[
r_1:= \frac{C_2^2 n^{1/4}}{4C_1}
\]
and $r_2:= C_2n^{1/4}$. Then by (\ref{fmu1}),
\[
f_\mu(r_1) \le\frac{C_2^2 n^{1/2}}{4}.
\]
On the other hand, since $r_2\le n^{1/4}$, therefore, by (\ref{fmu2}),
\[
f_\mu(r_2) \ge\frac{C_2^2n^{1/2}}{2}.
\]
Since $r_1 < r_2$, Proposition~\ref{simple} shows that $t_\mu\ge r_1$.
\end{pf*}

Finally, we give the proof of Proposition~\ref{adcor}, which is an
easy corollary of Theorem~\ref{admissibility}.

\begin{pf*}{Proof of Proposition~\ref{adcor}}
By Theorem~\ref{admissibility}, there exists $\mu_0\in K$ such that
$R_2 (\mu_0)\ge CR_1(\mu_0)$, where $C$ is a universal constant. Therefore,
\[
\sup_{\mu\in K} R_2(\mu) \ge R_2(
\mu_0) \ge C R_1(\mu_0)\ge C\inf
_{\mu\in K} R_1(\mu) = CL\sup_{\mu\in K}
R_1(\mu).
\]
This completes the proof of the proposition.
\end{pf*}

We now turn to the proofs of the theorems from Section~\ref
{applications}. The first goal is to prove Theorem~\ref{lassothm}. Let
us begin with some basic facts about Gaussian random variables.
%

\begin{lmm}[(Gaussian tails)]\label{gausstails}
Let $V$ be a standard Gaussian random variable. Then for any $x> 0$,
\begin{eqnarray*}
\biggl(\frac{1}{x}-\frac{1}{x^3}\biggr)\frac{2e^{-x^2/2}}{\sqrt{2\pi}}&\le& \pp\bigl(|V| > x\bigr) \le \frac{2e^{-x^2/2}}{x\sqrt{2\pi}}  , \\
\ee\bigl(|V|; |V|>x\bigr) &=& \frac{2e^{-x^2/2}}{\sqrt{2\pi}}\quad  \mbox{and} \\
\ee\bigl(V^2; |V|> x\bigr)& \le &\frac{2(x^2+1)e^{-x^2/2}}{x\sqrt{2\pi}}  .
\end{eqnarray*}
\end{lmm}

\begin{pf}
The upper bound in the first inequality is well known as the Mills ratio upper bound for the Gaussian tail. To prove this, just note that
\[
\pp\bigl(|V|> x\bigr) = 2\int_x^\infty \frac{e^{-y^2/2}}{\sqrt{2\pi}}   \,dy \le 2\int_x^\infty \frac{ye^{-y^2/2}}{x\sqrt{2\pi}}   \,dy = \frac{2e^{-x^2/2}}{x\sqrt{2\pi}}  .
\]
For the lower bound, we apply integration by parts two times to get
\[
\int_x^\infty e^{-y^2/2} \,dy = \biggl(\frac{1}{x}-\frac{1}{x^3}\biggr) e^{-x^2/2} + \int_x^\infty\frac{3e^{-y^2/2}}{y^4} \,dy  .
\]
For the second assertion, note that
\[
\ee\bigl(|V|; |V|>x\bigr) = 2\int_x^\infty \frac{ye^{-y^2/2}}{\sqrt{2\pi}} \,dy = \frac{2e^{-x^2/2}}{\sqrt{2\pi}}  .
\]
Finally, for the third claim, note that
\[
\ee\bigl(V^2; |V|> x\bigr) = 2\int_x^\infty \frac{y^2e^{-y^2/2}}{\sqrt{2\pi}}
= \frac{2x e^{-x^2/2}}{\sqrt{2\pi}} + 2\int_x^\infty  \frac{e^{-y^2/2}}{\sqrt{2\pi}}  \,dy
\]
and apply the first inequality to bound the second term on the right-hand side.
\end{pf}

%
\begin{lmm}[(Size of Gaussian maxima)]\label{maxsize}
Let $V_1,\ldots, V_n$ be standard Gaussian random variables, not
necessarily independent. Then
\[
\ee\Bigl(\max_{1\le i\le n}\llvert V_i\rrvert \Bigr)
\le\sqrt{2\log(2n)}.
\]
\end{lmm}

\begin{pf}
Take any $\beta> 0$. Then by Jensen's inequality,
\begin{eqnarray*}
\ee\Bigl(\max_{1\le i\le n}\llvert V_i\rrvert \Bigr) &=&
\frac{1}{\beta}\ee\bigl(\log e^{\beta
\max_{1\le i\le n}\llvert  V_i\rrvert  }\bigr)
\\
&\le&\frac{1}{\beta}\ee \Biggl(\log\sum_{i=1}^n
e^{\beta
\llvert  V_i\rrvert  } \Biggr) \le\frac{1}{\beta}\log\sum
_{i=1}^n \ee\bigl(e^{\beta
\llvert  V_i\rrvert  }\bigr)
\\
&\le&\frac{1}{\beta}\log\sum_{i=1}^n
\bigl(\ee\bigl(e^{\beta V_i}\bigr) + \ee \bigl(e^{-\beta V_i}\bigr)\bigr)=
\frac{\log(2n)}{\beta} + \frac{\beta}{2}.
\end{eqnarray*}
The proof is completed by taking $\beta= \sqrt{2\log(2n)}$.
\end{pf}
For any $n$ and $r$, let $C^r(\rr^n)$ be the set of $r$-times continuously
differentiable
functions from $\rr^n$ into $\rr$, and let $C^r_b(\rr^n)$ be the set
of all $g\in C^r(\rr^n)$ such that $g$ and all its derivatives up to
order $r$ are bounded. For any $g \in C^1_b(\rr)$, let $Ug$ be the
solution to the
differential equation
\[
f^\prime(x) - xf(x) = g(x) - \ee\bigl(g(V)\bigr),
\]
where $V\sim N(0,1)$. Explicitly, we have
\[
Ug(x) = e^{x^2/2}\int_{-\infty}^x
e^{-u^2/2}\bigl(g(u)-\ee\bigl(g(V)\bigr)\bigr)\,du.
\]
It is not difficult to prove that $Ug$ maps $C^1_b(\rr)$ into
$C^2_b(\rr)$. The following lemma is well known, and follows directly from
integration by parts:

%
\begin{lmm}\label{byparts}
Let $V = (V_1,\ldots,V_n)$ be a Gaussian random vector with zero
mean and arbitrary covariance matrix. Then for any $g \in C^1_b(\rr
^n)$ and any $i$, we have
\[
\ee\bigl(V_ig(V)\bigr) = \sum_{j=1}^n
\ee(V_iV_j) \ee \biggl(\frac{\partial
g}{\partial x_j}(V) \biggr).
\]
\end{lmm}
Using this, we easily get the following lemma.

%
\begin{lmm}\label{covlmm}
Take any $g_1,g_2\in C^2_b(\rr)$, and let $f_1 = Ug_1$, $f_2 =
Ug_2$. Suppose $V_1$ and $V_2$ are jointly Gaussian random variables with
$\ee(V_1) = \ee(V_2) = 0$, $\ee(V_1^2) = \ee(V_2^2) =1$ and
$\ee(V_1V_2) = \rho$. Then
\[
\cov\bigl(g_1(V_1),g_2(V_2)
\bigr) = \rho\ee\bigl(f_1(V_1)f_2(V_2)
\bigr) + \rho^2 \ee\bigl(f^\prime_1(V_1)f^\prime_2(V_2)
\bigr).
\]
\end{lmm}

\begin{pf}
Using Lemma~\ref{byparts} in two steps,
we have
\begin{eqnarray*}
\cov\bigl(g_1(V_1),g_2(V_2)
\bigr) &=& \ee \bigl(\bigl(f^\prime_1(V_1) -
V_1f_1(V_1)\bigr) \bigl(f^\prime_2(V_2)
- V_2f_2(V_2)\bigr) \bigr)
\\
&=& -\rho\ee \bigl(f_1(V_1) \bigl(f^{\prime\prime}_2(V_2)
- f_2(V_2) - V_2f^\prime_2(V_2)
\bigr) \bigr)
\\
&=& -\rho\ee \bigl(f_1(V_1) \bigl(f^{\prime\prime}_2(V_2)
- f_2(V_2)\bigr) \bigr)
\\
&&{}+ \rho\ee\bigl(f_1(V_1)f^{\prime\prime}_2(V_2)
\bigr) + \rho^2\ee\bigl(f^\prime_1(V_1)f^\prime_2(V_2)
\bigr)
\\
&=& \rho\ee\bigl(f_1(V_1)f_2(V_2)
\bigr) + \rho^2\ee\bigl(f^\prime_1(V_1)f^\prime_2(V_2)
\bigr).
\end{eqnarray*}
This completes the proof of the lemma.
\end{pf}
Using Lemma~\ref{covlmm}, we now prove the following set of
inequalities for additive functions of Gaussian random variables. This
is probably a new result.

%
\begin{lmm}\label{varlmm}
Let $V = (V_1,\ldots,V_n)$ be a Gaussian random vector with mean zero
and covariance matrix $\Sigma$. Let $\lambda_{\max}$ and $\lambda
_{\min}$ be the largest and smallest eigenvalues of $\Sigma$. Assume
that $\ee(V_i^2)=1$ for each $i$. Let $g_1,\ldots,g_n$ be
functions such that $\ee(g_i(V_i)^2) <\infty$
for each $i$. Then
\[
\lambda_{\min}\sum_{i=1}^n \var
\bigl(g_i(V_i)\bigr) \le\var \Biggl(\sum
_{i=1}^n g_i(V_i)
\Biggr) \le\lambda_{\max}\sum_{i=1}^n
\var\bigl(g_i(V_i)\bigr).
\]
\end{lmm}

\begin{pf}
First, let us make some
reductions. Recall that we have assumed that $\ee(V_i^2) =1 $
for each $i$. Next, note that if $g$ is a function such that
$\ee(g(Z)^2)<\infty$, where $Z\sim N(0,1)$, then there is a sequence
of step functions $\{g_n\}$ such that $g_n(Z)$ converges to $g(Z)$ in
$L^2$. Again, if $g$ is a step function, then there is a sequence
$\{g_n\}$ of
$C^1_b$ functions such that $g_n(Z)$ converges to $g(Z)$ in
$L^2$. Hence, assume without loss of generality that
$g_i$'s are elements of $C^1_b(\rr)$.

Now let $f_i:= Ug_i$ and $\sigma_{ij}:= \ee(V_iV_j)$. Let
$(Y_1,\ldots,Y_n)$ be an
independent copy of $(V_1,\ldots,V_n)$. Then by Lemma~\ref{covlmm}, we
have
\begin{eqnarray*}
\var \Biggl(\sum_{i=1}^n
g_i(V_i) \Biggr) &=& \sum
_{i,j} \bigl(\sigma_{ij}\ee\bigl(f_i(V_i)f_j(V_j)
\bigr) + \sigma_{ij}^2 \ee\bigl(f^\prime_i(V_i)f^\prime_i(V_j)
\bigr)\bigr)
\\
&=&\ee \biggl(\sum_{i,j}\sigma_{ij}
\bigl(f_i(V_i)f_j(V_j) +
Y_if^\prime_i(V_i)Y_jf^\prime_j(V_j)
\bigr) \biggr)
\\
&\le&\lambda_{\max} \ee \Biggl(\sum_{i=1}^n
\bigl(f_i(V_i)^2 + Y_i^2
f^\prime_i(V_i)^2 \bigr) \Biggr)
\\
&=& \lambda_{\max}\sum_{i=1}^n
\ee\bigl(f_i(V_i)^2 + f^\prime
_i(V_i)^2\bigr).
\end{eqnarray*}
But by Lemma~\ref{covlmm}, $\var(g_i(V_i)) = \ee(f(V_i)^2) +
\ee(f^\prime_i(V_i)^2)$. This gives the upper bound. The lower bound
follows similarly.
\end{pf}

We need a few more lemmas before proving Theorem~\ref{lassothm}. Let
all notation be as in the statement of the theorem. Additionally, let
$S:= \{i\dvtx \beta_i\ne0\}$, and let $V:= n^{-1/2}X^T Z$. Then $V$ is a
Gaussian random vector with mean zero and covariance matrix $\Sigma$.

%
\begin{lmm}\label{lassolmm1}
Suppose that $\delta> 0$. Take any $\alpha> 0$. Then there is a
constant $c_3$ depending only on $\alpha$, $\delta$, $a$, $b$, $s$,
$r$ and $L$ such that
\[
f_\mu\bigl(n^\alpha\bigr) \le c_3\sqrt{n}(
\log n)^{1/4} + c_3 n^\alpha\sqrt {\log n} + 2
\delta\sqrt{\alpha n\log n} - \frac{n^{2\alpha}}{2}.
\]
\end{lmm}

\begin{pf}
Throughout this proof, we will use $C_1, C_2,\ldots$ to denote
constants that may depend only on $\alpha$, $\delta$, $a$, $b$, $s$,
$r$ and $L$. Let $K_0'$ be the set of all $\gamma\in K_0$ such that $\llVert
X\gamma-X\beta\rrVert  \le n^\alpha$. Let
\[
M:= \sup_{\gamma\in K_0'} Z\cdot(X\gamma-X\beta) = \sqrt{n}\sup
_{\gamma\in K_0'} V\cdot(\gamma-\beta).
\]
Note that for any $\gamma\in K_0'$,
%
\begin{equation}
\label{betagamma} na\llVert \gamma- \beta\rrVert ^2\le\llVert X\gamma-X
\beta\rrVert ^2 \le n^{2\alpha}.
\end{equation}
Next, note that
%
\begin{eqnarray}
\label{vbeta} V\cdot(\gamma-\beta) &\le&\sum_{i\in S}
\llvert V_i\rrvert \llvert \gamma_i -
\beta_i\rrvert + \sum_{i\notin S} \llvert
V_i\rrvert \llvert \gamma_i\rrvert.
\end{eqnarray}
Now, by (\ref{betagamma}),
\begin{eqnarray*}
\sum_{i\in S} \llvert V_i\rrvert \llvert
\gamma_i-\beta_i\rrvert &\le& \biggl(\sum
_{i\in S} V_i^2 \sum
_{i\in S}(\gamma_i-\beta_i)^2
\biggr)^{1/2}
\\
&\le& \biggl(\sum_{i\in S} V_i^2
\biggr)^{1/2}\llVert \gamma-\beta\rrVert
\\
& \le& \frac{n^\alpha}{\sqrt{na}} \biggl(
\sum_{i\in S} V_i^2
\biggr)^{1/2}.
\end{eqnarray*}
Since $\ee(V_i^2)=1$ for each $i$, this shows that
%
\begin{eqnarray}
\label{vineq1} \ee \biggl(\sup_{\gamma\in K_0'} \sum
_{i\in S} \llvert V_i\rrvert \llvert \gamma
_i-\beta_i\rrvert \biggr) &\le& n^\alpha
\sqrt{\frac{s}{na}}.
\end{eqnarray}
Define the random set
\[
T:= \bigl\{i\notin S\dvtx  \llvert V_i\rrvert \ge2\sqrt{\alpha\log n}
\bigr\}.
\]
Then by (\ref{betagamma}) and the fact that $\llvert  \gamma\rrvert  _1\le L$,
\begin{eqnarray*}
\sum_{i\notin S} \llvert V_i\rrvert \llvert
\gamma_i\rrvert &\le&\sum_{i\in T} \llvert
V_i\rrvert \llvert \gamma_i\rrvert + 2\sqrt{\alpha
\log n}\sum_{i\notin S\cup T} \llvert \gamma_i
\rrvert
\\
&\le& \biggl(\sum_{i\in T} V_i^2
\biggr)^{1/2}\llVert \gamma-\beta\rrVert + 2\sqrt{\alpha\log n} \biggl(L
- \sum_{i\in S} \llvert \gamma_i\rrvert
\biggr)
\\
&\le& \biggl(\sum_{i\in T} V_i^2
\biggr)^{1/2}\frac{n^\alpha}{\sqrt
{na}} + 2\sqrt{\alpha\log n} \biggl(\delta+
\sum_{i\in S} \llvert \gamma _i-
\beta_i\rrvert \biggr).
\end{eqnarray*}
Again, by the Cauchy--Schwarz inequality and (\ref{betagamma}),
\begin{eqnarray*}
\sum_{i\in S} \llvert \gamma_i-
\beta_i\rrvert &\le&\sqrt{s}\llVert \gamma-\beta\rrVert \le
n^\alpha\sqrt{\frac{s}{na}}.
\end{eqnarray*}
From the last two displays, we get
\begin{eqnarray*}
\sum_{i\notin S} \llvert V_i\rrvert \llvert
\gamma_i\rrvert &\le& \biggl[ \biggl(\sum
_{i\in T} V_i^2 \biggr)^{1/2}
+ 2\sqrt{s\alpha\log n} \biggr]\frac{n^\alpha
}{\sqrt{na}} + 2\delta\sqrt{\alpha\log n}.
\end{eqnarray*}
Therefore, by Lemma~\ref{gausstails},
\begin{eqnarray*}
&&\ee \biggl(\sup_{\gamma\in K_0'} \sum_{i\notin S}
\llvert V_i\rrvert \llvert \gamma _i\rrvert \biggr)
\\
&&\qquad \le \Biggl[ \Biggl(\sum_{i=1}^p \ee
\bigl(V_i^2; \llvert V_i\rrvert \ge2\sqrt
{\alpha\log n} \bigr) \Biggr)^{1/2} + 2\sqrt{s\alpha\log n} \Biggr]
\frac{n^\alpha}{\sqrt{na}} + 2\delta\sqrt{\alpha\log n}
\\
&&\qquad \le \bigl[C_1n^{(1-2\alpha)/2}(\log n)^{1/4} + 2\sqrt{s
\alpha \log n} \bigr]\frac{n^\alpha}{\sqrt{na}} + 2\delta\sqrt{\alpha \log n}
\\
&&\qquad \le C_2(\log n)^{1/4} + C_3n^\alpha
\sqrt{\frac{\log n}{n}} + 2\delta\sqrt{\alpha\log n}.
\end{eqnarray*}
From the above display, and the inequalities (\ref{vbeta}) and (\ref
{vineq1}), we get
\begin{eqnarray*}
f_\mu\bigl(n^\alpha\bigr) &=& \ee(M) -\frac{n^{2\alpha}}{2} =
\sqrt{n} \ee \Bigl(\sup_{\gamma\in K_0'} V\cdot(\gamma-\beta) \Bigr) -
\frac
{n^{2\alpha}}{2}
\\
&\le& C_2\sqrt{n} (\log n)^{1/4} + C_3n^\alpha
\sqrt{\log n} +2\delta \sqrt{\alpha n \log n} - \frac{n^{2\alpha}}{2}.
\end{eqnarray*}
This completes the proof of the lemma.
\end{pf}

%
\begin{lmm}\label{lassolmm2}
Suppose that $\delta> 0$. Take any $0< \alpha_1<\alpha_2<1/2$. Then
there is a constant $c_4$ depending only on $\alpha_1$, $\alpha_2$,
$\delta$, $a$, $b$, $s$, $r$ and $L$ such that if $n> c_4$, then
\[
f_\mu\bigl(n^{\alpha_2}\bigr) \ge2\delta\sqrt{
\alpha_1n\log n} - 2n^{\alpha
_2} - \frac{n^{2\alpha_2}}{2}.
\]
\end{lmm}

\begin{pf}
Choose some $\alpha\in (\alpha_1, \alpha_2)$. Throughout this proof, we will use $C_1, C_2,\ldots$
to denote constants that may depend only on $\alpha, \alpha_1, \alpha_2, \delta, a, b, s, r$ and~$L$.

Let  $V$ and $T$ be as in the proof of Lemma \ref{lassolmm1}. Let  $K_0'$
and $M$ be as in the proof of Lemma \ref{lassolmm1}, with $\alpha$ replaced by $\alpha_2$.
Let us  make the following specific choice of $\gamma$:
\[
\gamma_i :=
\cases{
\operatorname{sign}(V_i)\delta/|T|,&\quad $\mbox{if } i\in T,$\vspace*{2pt}\cr
\beta_i, & \quad$\mbox{if } i \in S,$\vspace*{2pt}\cr
0, &\quad $\mbox{in all other cases.}$}
\]
Then note that
\begin{equation}\label{newineq0}
|\gamma|_1 \le |\beta|_1 + \delta = L  .
\end{equation}
(The above inequality is an equality if $T$ is nonempty, but we are allowing
for the possibility that $T$ may be empty.) Next, note that if $T$ is nonempty, then
\begin{equation}\label{newineq}
\|X\beta - X\gamma\|\le \sqrt{bn} \|\beta-\gamma\| \le \delta \sqrt{\frac{bn}{|T|}}
\end{equation}
and
\begin{equation}\label{vcdot}
V\cdot (\gamma-\beta) = \frac{\delta}{|T|} \sum_{i\in T} |V_i|  .
\end{equation}
By Lemma \ref{gausstails},
\begin{equation}\label{et1}
\biggl(1-\frac{1}{4\alpha \log n}\biggr)\frac{(p-s)n^{-2\alpha}}{\sqrt{2\pi\alpha\log n}}\le \ee|T| \le \frac{pn^{-2\alpha} }{\sqrt{2\pi\alpha\log n}}
\end{equation}
and
\begin{equation}\label{et2}
\ee\biggl( \sum_{i\in T} |V_i|\biggr) = \sum_{i\notin S} \ee\bigl(|V_i|; |V_i|\ge 2\sqrt{\alpha \log n}\bigr) = \frac{2(p-s)n^{-2\alpha}}{\sqrt{2\pi}}   .
\end{equation}
On the other hand, by Lemma \ref{varlmm},
\begin{equation}\label{et3}
\var\bigl(|T|\bigr) \le b\sum_{i\notin S} \pp\bigl(|V_i|\ge 2\sqrt{\alpha \log n}\bigr) \le \frac{C_1pn^{-2\alpha}}{\sqrt{\log n}}
\end{equation}
and
\begin{equation}\label{et4}
\var\biggl( \sum_{i\in T} |V_i|\biggr) = b \sum_{i\notin S} \ee\bigl(V_i^2; |V_i|\ge 2\sqrt{\alpha \log n}\bigr) \le C_2pn^{-2\alpha}\sqrt{\log n}  .
\end{equation}
Let $\ep'$ be a positive constant depending only on $\alpha$, $\alpha_1$, $\alpha_2$, $\delta$, $a$, $b$, $s$, $r$ and $L$.
The value of $\ep'$ will be determined later. As a consequence of  \eqref{newineq}, \eqref{et1}, \eqref{et2}, \eqref{et3},
\eqref{et4}, the fact that $\alpha_2 <1/2$, and Chebychev's inequality, it follows that there exists $C_3$ depending only
on $\alpha$, $\alpha_1$, $\alpha_2$,  $\delta$, $a$, $b$, $s$, $r$ and $L$ and our choice of $\ep'$, such that if $n> C_3$, then
\begin{eqnarray*}
&&\pp\biggl(\bigl(1-\ep'\bigr) \biggl(1-\frac{1}{4\alpha \log n}\biggr)\frac{(p-s)n^{-2\alpha}}{\sqrt{2\pi\alpha\log n}}\le
 |T|\le \bigl(1+ \ep'\bigr)\frac{pn^{-2\alpha} }{\sqrt{2\pi\alpha\log n}}\\
&& \hspace*{138pt}\mbox{and } \sum_{i\in T} |V_i|\ge \bigl(1-\ep'^2\bigr)\frac{2(p-s)n^{-2\alpha}}{\sqrt{2\pi}}\biggr) \ge \frac{1}{2}.
\end{eqnarray*}
Note that if $C_3$ is chosen large enough, and $|T|$ indeed turns out to be bigger than the lower
bound on $|T|$ in the above expression, then $|T|\ge \delta^2 b n^{1-2\alpha_2}$ since $\alpha_2 > \alpha$.
Thus, under this circumstance \eqref{newineq0} and \eqref{newineq} imply that $\gamma\in K_0'$.
Combined with  \eqref{vcdot} and the lower bound on the probability displayed above, this gives
\[
\pp\biggl(M \ge \frac{(1-\ep') (p-s)}{p}2\delta \sqrt{\alpha n\log n}\biggr)\ge \frac{1}{2}  .
\]
By the concentration of Gaussian maxima (Lemma \ref{gaussmax}) and the above inequality, it follows that
\[
\ee(M)\ge \frac{(1-\ep') (p-s)}{p}2\delta \sqrt{\alpha n\log n} - 2n^{\alpha_2}  .
\]
The proof is now completed by taking $\ep'$ small enough and $C_3$ large enough to satisfy the required inequality.
\end{pf}

%
\begin{lmm}\label{lassolmm3}
Suppose that $\delta= 0$. Then there is a constant $c_5$ depending
only on $a$, $b$, $s$, $r$ and $L$ such that for any $u > 0$,
\[
f_\mu (u\sqrt{\log n} )\le c_5 u \log n -
\frac{u^2\log
n}{2}.
\]
\end{lmm}

\begin{pf}
Throughout this proof, we will use $C_1, C_2,\ldots$ to denote
constants that may depend only on $\delta$, $a$, $b$, $s$, $r$ and
$L$. Fix $u>0$ and let $K_0'$ be the set of all $\gamma\in K_0$ such
that $\llVert  X\gamma-X\beta\rrVert  \le u\sqrt{\log n}$. Let $M$ and $V$ be as
in the proof of Lemma~\ref{lassolmm1}. Additionally, let $G:= \max_{1\le i\le p} \llvert  V_i\rrvert  $.

Take any $\gamma\in K_0'$. Note that the inequality (\ref{vbeta})
from the proof of Lemma~\ref{lassolmm1} is still valid, and that (\ref
{betagamma}) and (\ref{vineq1}) are also valid, after replacing
$n^\alpha$ with $u\sqrt{\log n}$. In addition to that, note that by
Lemma~\ref{maxsize},
\begin{eqnarray*}
\ee \biggl(\sum_{i\notin S} \llvert V_i
\rrvert \llvert \gamma_i\rrvert \biggr) &\le&\ee(G) \sum
_{i\notin S} \llvert \gamma_i\rrvert
\\
&\le&\sqrt{2\log(2p)} \biggl(L-\sum_{i\in S} \llvert
\gamma_i\rrvert \biggr)
\\
&=& \sqrt{2\log(2p)} \sum_{i\in S} \bigl(\llvert
\beta_i\rrvert - \llvert \gamma_i\rrvert \bigr)\le
\sqrt{2\log(2p)} \sum_{i\in S} \llvert
\beta_i - \gamma_i\rrvert
\\
&\le&\sqrt{2s\log(2p)} \llVert \beta-\gamma\rrVert \le\frac{C_1 u\log
n}{\sqrt{n}}.
\end{eqnarray*}
Combining the above observations, we get
\[
\ee(M)\le C_2u\sqrt{\log n} + C_1 u\log n -
\frac{u^2\log n}{2}.
\]
This completes the proof of the lemma.
\end{pf}

\begin{lmm}\label{lassolmm4}
Suppose that $\delta< 0$. Then there are positive constants $c_6$ and
$c_7$ depending only on $\delta$, $a$, $b$, $s$, $r$ and $L$ such that
$c_6\sqrt{n}\le t_\mu\le c_7 \sqrt{n}$.
\end{lmm}

\begin{pf}
Note that for any $\gamma\in K_0$,
\begin{eqnarray*}
\llVert X\gamma-X\beta\rrVert ^2 &\ge& na\llVert \gamma-\beta\rrVert
^2 \ge na\sum_{i\in S} (
\gamma_i-\beta_i)^2
\\
&\ge&\frac{na}{s} \biggl(\sum_{i\in S} \llvert
\gamma_i-\beta_i\rrvert \biggr)^2\ge
\frac{na\delta^2}{s}.
\end{eqnarray*}
This shows that there is a small enough $C_1$ depending only on $\delta
$, $a$ and $s$ such that $f_\mu(t)=-\infty$ if $t\le C_1\sqrt{n}$.
By Proposition~\ref{simple} and the fact that $f_\mu(t)$ is finite
for at least one $t$ (from Theorem~\ref{mainest}), this implies the
lower bound on $t_\mu$.

Next, note that since $0\in K_0$,
\begin{eqnarray*}
\llVert \mu-\hmm\rrVert &\le&\llVert \mu-Y\rrVert + \llVert Y-\hmm\rrVert
\\
&\le&\llVert \mu-Y\rrVert + \llVert Y\rrVert
\\
&\le&2\llVert \mu-Y\rrVert + \llVert \mu\rrVert.
\end{eqnarray*}
But $\ee\llVert  \mu-Y\rrVert  ^2 = n$ and
\[
\llVert \mu\rrVert =\llVert X\beta\rrVert \le\sqrt{nb} \llVert \beta\rrVert \le
\sqrt{nb} \llvert \beta\rrvert _1\le\sqrt{nb} L.
\]
Thus, $\ee\llVert  \mu-\hmm\rrVert  ^2 \le(8 + 2bL^2) n$. By Corollary~\ref
{maincor}, this shows that $t_\mu\le C_2\sqrt{n}$ for some constant
$C_2$ depending only on $b$ and $L$.
This completes the proof of the lemma.
\end{pf}
We are now ready to prove Theorem~\ref{lassothm}.

\begin{pf*}{Proof of Theorem~\ref{lassothm}}
First, suppose that $\delta> 0$. Take any $0< \alpha< \alpha_1
<\alpha_2< 1/4$. By Lemmas~\ref{lassolmm1}~and~\ref{lassolmm2},
it follows that if $n$ is large enough (depending only on $\alpha$,
$\alpha_1$, $\alpha_2$, $\delta$, $a$, $b$, $s$, $r$ and $L$), then
$f_\mu(n^\alpha) \le f_\mu(n^{\alpha_2})$ and, therefore, by
Proposition~\ref{simple}, $t_\mu\ge n^\alpha$. Next, take any
$\alpha> 1/4$. Lemma~\ref{lassolmm1} implies that if $n$ is large
enough, then $f_\mu(n^\alpha)\le0$ and, therefore, by Proposition
\ref{simple}, $t\le n^\alpha$.

If $\delta= 0$, the conclusion follows directly from a combination of
Lemma~\ref{lassolmm3} and Proposition~\ref{simple}. If $\delta< 0$,
simply invoke Lemma~\ref{lassolmm4}.
\end{pf*}

Our final task is to prove Theorem~\ref{convthm}. As before, we need
some standard results and notations from the literature.

If $\mf$ is a subset of a normed space with norm $\llVert  \cdot\rrVert  $ and
$\varepsilon$ is a positive real number, the \emph{covering number}
$N(\varepsilon, \mf, \llVert  \cdot\rrVert  )$ is defined as the minimum number of
open balls of radius $\varepsilon$ (with respect to the norm $\llVert  \cdot\rrVert
$) with centers in $\mf$ that are needed to cover~$\mf$.

The following result, known as ``Dudley's entropy bound,'' connects the
covering numbers of $\mf$ with the expected maximum of a certain
Gaussian process.

%
\begin{lmm}[(Dudley's entropy bound \cite{dudley})]\label{dudleylmm}
Let $\mf$ be as above. Suppose that $(X_f)_{f\in\mf}$ is a Gaussian
process on $\mf$ such that $\ee(X_f)=0$ for each $f\in\mf$, and
$\ee(X_f-X_g)^2=\llVert  f-g\rrVert  ^2$ for each $f,g\in\mf$. Then
\[
\ee \Bigl(\sup_{f\in\mf} X_f \Bigr) \le C \int
_0^{\operatorname{diam}(\mf)} \sqrt{\log N\bigl(\varepsilon, \mf, \llVert
\cdot\rrVert \bigr) } \,d\varepsilon,
\]
where $C$ is a universal constant.
\end{lmm}
Suppose now that $\mf$ is a set of functions from some set $S$ into
$\rr$, and $\llVert  \cdot\rrVert  $ is a norm on a vector space of functions
containing $\mf$. Suppose that $l$ and $u$ are two elements of $\mf$
such that $l\le u$ everywhere on $S$. If $\llVert  l-u\rrVert  \le\varepsilon$, then
the set of all $f\in\mf$ such that $l\le f\le u$ everywhere on $S$ is
called an $\varepsilon$-bracket, and is denoted by $[l,u]$. The \emph{bracketing number} $N_{[\,]}(\varepsilon, \mf, \llVert  \cdot\rrVert  )$ is the
minimum number of $\varepsilon$-brackets needed to cover $\mf$. It is
quite easy to see that
%
\begin{equation}
\label{brac} N\bigl(\varepsilon, \mf, \llVert \cdot\rrVert \bigr) \le
N_{[\,]}\bigl(2\varepsilon, \mf, \llVert \cdot\rrVert \bigr).
\end{equation}
The following result is quoted from van der Vaart and Wellner \cite{vdvwellner}, Theorem 2.7.5, page~159.

%
\begin{lmm}[(van der Vaart and Wellner \cite{vdvwellner})]\label{vdvlmm}
Let $P$ be any probability measure on $\rr$ and let $\llVert  \cdot\rrVert  _r$
denote the $L^r(P)$ norm. Let $\mf$ be the set of all monotone
functions from $\rr$ into $[0,1]$. Then for any $\varepsilon>0$,
\[
\log N_{[\,]}\bigl(\varepsilon, \mf, \llVert \cdot\rrVert _r
\bigr)\le C\varepsilon^{-1},
\]
where $C$ is a constant that depends on $r$ only.
\end{lmm}
The statement of Lemma~\ref{vdvlmm} has to be modified in a certain
way to suit our purpose in the proof of Theorem~\ref{convthm}. The
following lemma gives the modified statement.

%
\begin{lmm}\label{vdv2}
Take any two real numbers $a<b$, and a positive integer $n$. Let $Q$
denote the set of all vectors $\mu\in\rr^n$ such that
\[
a\le\mu_1\le\mu_2\le\cdots\le\mu_n\le b.
\]
Let $\llVert  \cdot\rrVert  $ denote the Euclidean norm on $Q$. Then for any $t>0$,
\[
\log N\bigl(t, Q, \llVert \cdot\rrVert \bigr) \le\frac{C\sqrt{n}(b-a)}{t},
\]
where $C$ is a universal constant.
\end{lmm}

\begin{pf}
First, assume that $a=0$ and $b=1$. Let
\[
\varepsilon:= \frac{t}{2\sqrt{n}}.
\]
Let $P$ be the uniform probability distribution on $[0,1]$, and let $\llVert
\cdot\rrVert  _{L^2(P)}$ denote the $L^2$ norm induced by $P$. Let $\mf$ be
the set of all monotone functions from $\rr$ into $[0,1]$. Let $\mg$
be a finite subset of $\mf$ such that for any $f\in\mf$ there exists
$g\in\mg$ such that $\llVert  f-g\rrVert  _{L^2(P)} \le\varepsilon$. By Lemma~\ref
{vdvlmm} and the inequality (\ref{brac}), $\mg$ can be chosen such
that $\log\llvert  \mg\rrvert  \le C\varepsilon^{-1}$, where $C$ is a universal constant.

Now take any $\mu\in Q$. Define a function $f^\mu\dvtx \rr\rightarrow
[0,1]$ as
\[
f^\mu(x) = \cases{ 0, &\quad if $x < 0$,
\cr
\mu_i, &
\quad if $(i-1)/n \le x< i/n$,
\cr
1, &\quad if $x \ge1$.}
\]
Then clearly $f^\mu\in\mf$. For each $g\in\mg$, inspect whether
there exists some $\mu\in Q$ such that $\llVert  f^\mu- g\rrVert  _{L^2(P)}<
\varepsilon$. If there exists such a $\mu$, choose one according to some
pre-specified rule and call it $\mu(g)$. Let $Q'$ be the subset of $Q$
consisting of all such $\mu(g)$. Then clearly $\llvert  Q'\rrvert  \le\llvert  \mg\rrvert  $. On the
other hand, for any $\mu\in Q$, there exists $g\in\mg$ such that $\llVert
f^\mu-g\rrVert  _{L^2(P)}< \varepsilon$. Consequently,
\[
\bigl\llVert f^\mu-f^{\mu(g)}\bigr\rrVert _{L^2(P)}< 2
\varepsilon.
\]
But
\begin{eqnarray*}
\bigl\llVert f^\mu-f^{\mu(g)}\bigr\rrVert _{L^2(P)}^2
&=& \int_0^1 \bigl(f^\mu(x)-f^{\mu
(g)}(x)
\bigr)^2 \,dx= \frac{1}{n}\sum_{i=1}^n
\bigl(\mu_i-\mu_i(g)\bigr)^2.
\end{eqnarray*}
Thus, $\llVert  \mu-\mu(g)\rrVert  = \sqrt{n}\llVert  f^\mu-f^{\mu(g)}\rrVert  _{L^2(P)} <
2\sqrt{n}\varepsilon= t$. This completes the proof of the lemma when
$a=0$ and $b=1$.

For general $a$ and $b$, let $l$ be the unique linear map that takes
$a$ to $0$ and $b$ to~$1$. Let $L\dvtx \rr^n \rightarrow\rr^n$ be the map
that applies $l$ to each coordinate. Given $t>0$, we now know that
there exists a set $Q_0\subseteq L(Q)$ of size $\le C\sqrt{n}(b-a)/t$
such that for any $\mu\in Q$, there exists $\nu\in Q_0$ satisfying
\[
\bigl\llVert L(\mu)-\nu\bigr\rrVert \le\frac{t}{b-a}.
\]
To complete the proof, note that $L^{-1}(Q_0)\subseteq Q$, and $\llVert  \mu
-L^{-1}(\nu)\rrVert  \le t$.
\end{pf}
We are now ready to prove Theorem~\ref{convthm}.

\begin{pf*}{Proof of Theorem~\ref{convthm}}
Fix $\mu\in K$. Let $l$ be a positive integer, to be chosen later. Let
$K'$ be the subset of $K$ consisting of all $\nu$ such that
\[
\nu_1\ge\mu_1 - 2^{l}, \nu_n \le
\mu_n + 2^l.
\]
Fix $t>0$. Let
\[
K'':= \bigl\{\nu\in K'\dvtx  \llVert \nu-
\mu\rrVert \le t\bigr\}
\]
and
\[
m:= \ee \Bigl(\sup_{\nu\in K''} Z\cdot(\nu-\mu) \Bigr).
\]
Given any $s> 0$, Lemma~\ref{vdv2} implies that there exists a set
$A\subseteq K'$ of size $\le\exp(C_02^lD\sqrt{n}/s)$ such that for
any $\nu\in K'$ there exists $\gamma\in A$ satisfying \mbox{$\llVert  \nu-\gamma
\rrVert  < s$}. Combined with Dudley's entropy bound (Lemma~\ref{dudleylmm}),
this gives
%
\begin{equation}
\label{isoineq0} m\le C_1\sqrt{2^lD} n^{1/4}\int
_0^t \frac{ds}{\sqrt{s}} = 2C_1
\sqrt{ 2^lDt} n^{1/4}.
\end{equation}
Now take any $\nu\in K$ such that $\llVert  \nu-\mu\rrVert  \le t$. For any $L > 0$,
\[
\bigl\llvert \bigl\{i\dvtx  \llvert \nu_i-\mu_i\rrvert >
L\bigr\}\bigr\rrvert \le\frac{1}{L^2}\sum_{i=1}^n
(\nu_i-\mu _i)^2 \le\frac{t^2}{L^2}.
\]
Consequently, if $r(L)$ is the largest $i$ such that $\llvert  \nu_i-\mu
_i\rrvert  \le L$, then
\[
r(L)\ge n - \frac{t^2}{L^2}.
\]
Similarly, if $s(L)$ is the smallest $i$ such that $\llvert  \nu_i-\mu_i\rrvert   \le
L$, then
\[
s(L)\le1 + \frac{t^2}{L^2}.
\]
Define $\nu'$ as
\[
\nu'_i:= \cases{ \mu_i +
2^l, &\quad if $ i> r\bigl(2^l\bigr)$,
\vspace*{3pt}\cr
\mu_i - 2^l, &\quad if $i< s\bigl(2^l
\bigr)$,
\vspace*{3pt}\cr
\nu_i, &\quad if $s\bigl(2^l\bigr) \le i
\le r\bigl(2^l\bigr)$.}
\]
Since $\nu_{r(2^l)} \le\mu_{r(2^l)}+2^l$ and $\nu_{s(2^l)}\ge\mu
_{s(2^l)} -2^l$, we see that $\nu'\in K$. Again by definition it is
clear that $\nu_n'\le\mu_n+2^l$ and $\nu_1' \ge\mu_1-2^l$.
Therefore, $\nu'\in K'$. Finally, note that for any $i$, $\llvert  \mu_i-\nu
'_i\rrvert  \le\llvert  \mu_i-\nu_i\rrvert  $, implying that $\nu'\in K''$. Thus,
%
\begin{equation}
\label{isoineq1} Z\cdot\bigl(\nu'-\mu\bigr)\le\sup
_{\gamma\in K''} Z\cdot(\gamma-\mu).
\end{equation}
Next, note that
\begin{eqnarray*}
Z\cdot\bigl(\nu-\nu'\bigr) &\le&\sum_{i>r(2^l)}
\llvert Z_i\rrvert \bigl\llvert \nu_i-
\nu_i'\bigr\rrvert + \sum
_{i< s(2^l)} \llvert Z_i\rrvert \bigl\llvert
\nu_i-\nu_i'\bigr\rrvert
\\
&\le&\sum_{k=l}^\infty\sum
_{r(2^k)< i\le r(2^{k+1})} \llvert Z_i\rrvert \bigl\llvert \nu
_i-\nu_i'\bigr\rrvert +\sum
_{k=l}^\infty\sum_{s(2^{k+1})\le i <s(2^{k})}
\llvert Z_i\rrvert \bigl\llvert \nu_i-
\nu_i'\bigr\rrvert
\\
&\le&\sum_{k=l}^\infty\sum
_{r(2^k)< i\le r(2^{k+1})} \llvert Z_i\rrvert 2^{k+2}+\sum
_{k=l}^\infty\sum
_{s(2^{k+1})\le i < s(2^{k})} \llvert Z_i\rrvert 2^{k+2}
\\
&\le&\sum_{k=l}^\infty\sum
_{i> n - t^2/2^{2k}} \llvert Z_i\rrvert 2^{k+2}+\sum
_{k=l}^\infty\sum
_{i < 1+ t^2/2^{2k}} \llvert Z_i\rrvert 2^{k+2}.
\end{eqnarray*}
This shows that
%
\begin{eqnarray}
\label{isoineq2} \ee \Bigl(\sup_{\nu\in K\dvtx    \llVert  \nu-\mu\rrVert  \le t} Z\cdot\bigl(\nu -
\nu'\bigr) \Bigr) &\le&\sum_{k=l}^\infty
\frac{C_2t^2}{2^k}= \frac{C_2
t^2}{2^{l-1}}.
\end{eqnarray}
Combining (\ref{isoineq0}), (\ref{isoineq1}) and (\ref{isoineq2}) gives
\begin{eqnarray*}
\ee \Bigl(\sup_{\nu\in K\dvtx    \llVert  \nu-\mu\rrVert  \le t} Z\cdot(\nu -\mu) \Bigr) &\le&\ee \Bigl(
\sup_{\nu\in K\dvtx    \llVert  \nu-\mu\rrVert
\le t} Z\cdot\bigl(\nu'-\mu\bigr) \Bigr)
\\
&&{}+ \ee \Bigl(\sup_{\nu\in K\dvtx    \llVert  \nu-\mu\rrVert  \le t} Z\cdot\bigl(\nu-
\nu'\bigr) \Bigr)
\\
&\le&\ee \Bigl(\sup_{\gamma\in K''} Z\cdot(\gamma-\mu) \Bigr) +
\frac{C_2t^2}{2^{l-1}}
\\
&\le&2C_1\sqrt{2^lDt} n^{1/4} +
\frac{C_2t^2}{2^{l-1}}.
\end{eqnarray*}
Now choose $l$ so large that $C_22^{-(l-1)}\le1/4$. With this choice
of $l$, the above inequality implies that
%
\begin{equation}
\label{fmuup} f_\mu(t)\le C_3 \sqrt{Dt}
n^{1/4}-\frac{t^2}{4}.
\end{equation}
In particular, $f_\mu(r) \le0$, where $r = (4C_3 \sqrt{D}
n^{1/4})^{2/3}$. By Proposition~\ref{simple}, this implies that $t_\mu
\le r$. This completes the proof of the upper bound for $t_\mu$ in the
statement of the theorem.

Next, fix $t\in[Bn^{-1/2}, \sqrt{n}]$. Let $k:= \lceil t\sqrt
{n}/B\rceil$ and $m:= \lfloor n/k\rfloor$. For $j=1,2,\ldots,m$, let
\[
S_j:= \sum_{(j-1)k< i\le jk} Z_i,\qquad
a_j:= \mu_{(j-1)k+1},\qquad b_j:= \mu_{jk}
\]
and if $mk< n$, let
\[
S_{m+1}:= \sum_{mk< i\le n} Z_i,\qquad
a_{m+1}:= \mu_{mk+1},\qquad b_{m+1}:= \mu_n.
\]
For each $i$, let
\[
\nu_i:= \frac{a_j+b_j}{2} \qquad\mbox{ if } (j-1)k < i \le jk.
\]
Additionally, define
\[
\gamma_i:= \cases{ a_j, &\quad if $(j-1)k < i \le
jk$ and $S_j < 0$,
\vspace*{2pt}\cr
b_j, &\quad if $(j-1)k < i \le
jk$ and $S_j > 0$.}
\]
Notice that for each $i$,
\[
\llvert \gamma_i-\mu_i\rrvert \le\frac{Bk}{n}
\le\frac{t}{\sqrt{n}}.
\]
Consequently,
\[
\llVert \gamma- \mu\rrVert \le t.
\]
Moreover, $\gamma\in K$. Next, note that
\begin{eqnarray*}
Z\cdot(\gamma-\nu) &=& \frac{1}{2}\sum_{j=1}^{m+1}
\llvert S_j\rrvert (b_j-a_j)
\ge\frac{Ak}{2n}\sum_{j=1}^{m}
\llvert S_j\rrvert.
\end{eqnarray*}
Therefore,
\begin{eqnarray*}
\ee \Bigl(\sup_{\theta\in K\dvtx   \llVert  \theta-\mu\rrVert  \le t} Z\cdot (\theta-\mu) \Bigr) &=& \ee \Bigl(
\sup_{\theta\in K\dvtx   \llVert
\theta-\mu\rrVert  \le t} Z\cdot(\theta-\nu) \Bigr)
\\
&\ge&\ee\bigl(Z\cdot(\gamma-\nu)\bigr)
\\
&\ge&\frac{Ak}{2n}\sum_{j=1}^m \ee
\llvert S_j\rrvert \ge\frac{C_4Akm\sqrt
{k}}{n}
\\
&\ge& C_5A\sqrt{k}\ge C_5A B^{-1/2}
t^{1/2} n^{1/4}.
\end{eqnarray*}
Thus,
%
\begin{equation}
\label{fmudown} f_\mu(t)\ge C_5 A B^{-1/2}
t^{1/2}n^{1/4}-\frac{t^2}{2}.
\end{equation}
Let $\alpha$ and $\beta$ be two positive constants, to be chosen
later. Let
\[
r_1:= \alpha A^{8/3} B^{-4/3}D^{-1}
n^{1/6},\qquad r_2:= \beta A^{2/3} B^{-1/3}
n^{1/6}.
\]
Then by (\ref{fmudown}),
\[
f_\mu(r_2)\ge\bigl(C_5\sqrt{\beta}-
\beta^2/2\bigr) A^{4/3} B^{-2/3} n^{1/3},
\]
and by (\ref{fmuup}),
\[
f_\mu(r_1)\le C_3 \alpha^{1/2}
A^{4/3} B^{-2/3} n^{1/3}.
\]
Suppose that $A> 0$. Choosing $\beta$ sufficiently small, and then
choosing $\alpha$ even smaller (depending on $\beta$), it is now easy
to arrange that $r_1 < r_2$ and $f_\mu(r_1)\le f_\mu(r_2)$. By
Proposition~\ref{simple}, this implies that $t_\mu\ge r_1$. If $A=0$,
the lower bound in the statement of the theorem is automatically true.
\end{pf*}


\section*{Acknowledgments}
The author thanks Bodhisattva Sen for introducing him to this area and many useful discussions,
Xi Chen for pointing out a mistake in the proof of Lemma \ref{lassolmm2} in an earlier draft,
Joel Tropp for pointing out the relevant signal processing literature, and Adityanand Guntuboyina
and Sara van de Geer for helpful comments. The author also thanks the anonymous
referees and the  Associate Editor for several useful suggestions.


%

\printaddresses
\end{document}